\documentclass[11pt,letterpaper]{article}
\usepackage{amsfonts,amssymb,amsmath,amsthm,times1}
\usepackage{graphicx}
\makeatletter

\renewcommand{\@listI}{}
\renewcommand{\@listI}{%
\setlength{\labelwidth}{225pc}
\setlength{\leftmargin}{15pt}
}

\renewcommand{\@biblabel}[1]{#1.}

\makeatother

\newtheorem{theorem}{Theorem}

\theoremstyle{definition}

\DeclareMathOperator{\EW}{\mathit{E}}
\DeclareMathOperator{\PW}{\mathit{P}}
\DeclareMathOperator{\EE}{\mathbb{E}}
\DeclareMathOperator{\PP}{\mathbb{P}}
\DeclareMathOperator{\QQ}{\mathbb{Q}}

\DeclareMathOperator{\Prob}{\mbox{\boldmath$P$}}
\DeclareMathOperator{\Exp}{\mbox{\boldmath$E$}}

\newcommand{\e}{{\mathrm e}}
\newcommand{\dif}{{\mathrm d}}

\newcommand{\RR}{\mathbb{R}}
\newcommand{\ZZ}{\mathbb{Z}}

\newcommand{\calF}{\mathcal{F}}

\begin{document}
\noindent \twocolumn[{\small {\it Encyclopedia of Mathematical
Physics} \,(J.-P. Fran\c{c}oise, G. Naber, and S.~T. Tsou, eds.)\\
Vol. 4, pp. 353--371. Elsevier, Oxford, 2006.

\vspace{2pc}

\sffamily\bfseries\LARGE Random Walks in Random Environments}

\vspace{1pc}]

\noindent {\sffamily\bfseries\small L. V. Bogachev,
\normalfont\sffamily University of Leeds, Leeds, UK} 

\vspace{.3pc} \noindent {\sffamily\scriptsize\copyright\ 2006
Elsevier Ltd. All rights reserved.}

\renewcommand{\thefootnote}{}


\section{Introduction}\label{sec:Intro}

Random walks provide a simple conventional model to describe various
transport processes, for example propagation of heat or diffusion of
matter through a medium (for a general reference see, e.g., Hughes
(1995)). However, in many practical cases the medium where the
system evolves is highly irregular, due to factors such as defects,
impurities, fluctuations etc. It is natural to model such
irregularities as \emph{random environment}, treating the observable
sample as a statistical realization of an ensemble, obtained by
choosing the local characteristics of the motion (e.g., transport
coefficients and driving fields) at random, according to a certain
probability distribution.

In the random walks context, such models are referred to as
\emph{Random Walks in Random Environments} (RWRE). This is a
relatively new chapter in applied probability and physics of
disordered systems initiated in the 1970s. Early interest in RWRE
models was motivated by some problems in biology, crystallography
and metal physics, but later applications have spread through
numerous areas (see review papers by Alexander \emph{et al.} (1981),
Bou\-chaud and Georges (1990), and a comprehensive monograph by
Hughes (1996)). After 30 years of extensive work, RWRE remain a very
active area of research, which has been a rich source of hard and
challenging questions and has already led to many surprising
discoveries, such as subdiffusive behavior, trapping effects,
localization, etc. It is fair to say that the RWRE paradigm has
become firmly established in physics of random media, and its
models, ideas, methods, results, and general effects have become an
indispensable part of the standard tool kit of a mathematical
physicist.

One of the central problems in random media theory is to establish
conditions ensuring homogenization, whereby a given stochastic
system evolving in a random medium can be adequately described, on
some spatial-temporal scale, using a suitable effective system in a
homogeneous (non-random) medium. In particular, such systems would
exhibit classical diffusive behavior with effective drift and
diffusion coefficient. Such an approximation, called \emph{effective
medium approximation} (EMA), may be expected to be successful for
systems exposed to a relatively small disorder of the environment.
However, in certain circumstances EMA may fail due to atypical
environment configurations (``large deviations'') leading to various
anomalous effects. For instance, with small but positive probability
a realization of the environment may create ``traps'' that would
hold the particle for anomalously long time, resulting in the
\emph{sub\-dif\-fusive} behavior, with the mean square displacement
growing slower than linearly in time.

RWRE models have been studied by various non-rigorous methods
including Monte Carlo simulations, series expansions, and the
renormalization group techniques (see more details in the above
references), but only few models have been analyzed rigorously,
especially in dimensions greater than one. The situation is much
more satisfactory in the one-dimensional case, where the
mathematical theory has matured and the RWRE dynamics has been
understood fairly well.

The goal of this article is to give a brief introduction to the
beautiful area of RWRE. The principal model to be discussed is a
random walk with nearest-neighbor jumps in independent identically
distributed (i.i.d.) random environment in one dimension, although
we shall also comment on some generalizations. The focus is on
rigorous results; however, heuristics will be used freely to
motivate the ideas and explain the approaches and proofs. In a few
cases, sketches of the proofs have been included, which should
help appreciate the flavor of the results and methods.

\subsection{Ordinary Random Walks: A Reminder}\label{sec:1.1}

To put our exposition in perspective, let us give a brief account
of a few basic concepts and facts for \emph{ordinary random
walks}, that is, evolving in a non-random environment (see further
details in Hughes 1995). In such models, space is modelled using
a suitable graph, e.g., a $d$-dimensional integer lattice $\ZZ^d$,
while time may be discrete or continuous. The latter distinction
is not essential, and in this article we will mostly focus on the
discrete-time case. The random mechanism of spatial motion is then
determined by the given transition probabilities (probabilities of
jumps) at each site of the graph. In the lattice case, it is
usually assumed that the walk is translation invariant, so that at
each step distribution of jumps is the same, with no regard to the
current location of the walk.

In one dimension ($d=1$), the \emph{simple (nearest-neighbor) random
walk} may move one step to the right or to the left at a time, with
some probabilities $p$ and $q=1-p$, respectively. An important
assumption is that only the current location of the walk determines
the random motion mechanism, whereas the past history is not
relevant. In terms of probability theory, such a process is referred
to as \emph{Markov chain}. Thus, assuming that the walk starts at
the origin, its position after $n$ steps can be represented as the
sum of consecutive displacements, $X_n=Z_1+\dots+Z_n$, where $Z_i$
are independent random variables with the same distribution
$P\{Z_i=1\}=p$, \,${P\{Z_i=-1\}=q}$.

The strong law of large numbers (LLN) states that almost surely
(i.e., with probability $1$)
\begin{equation}\label{LLN}
\lim_{n\to\infty}\frac{X_n}{n}=\EW Z_1=p-q,\quad \PW\text{-a.s.}
\end{equation}
where $\EW$ denotes expectation (mean value) with respect to
$\PW$. This result shows that the random walk moves with the
asymptotic average velocity close to $p-q$. It follows that if
$p-q\ne0$ then the process $X_n$, with probability 1, will
ultimately drift to infinity (more precisely, $+\infty$ if $p-q>0$
and $-\infty$ if $p-q<0$). In particular, in this case the random
walk may return to the origin (and in fact visit \emph{any} site
on $\ZZ$) only finitely many times. Such behavior is called
\emph{transient}. However, in the symmetric case (i.e., $p=q=0.5$)
the average velocity vanishes, so the above argument fails. In
this case the walk behavior appears to be more complicated, as it
makes increasingly large excursions both to the right and to the
left, so that $\varlimsup_{n\to\infty} X_n=+\infty$,
\,$\varliminf_{n\to\infty} X_n=-\infty$ ($\PW$-a.s.). This implies
that a symmetric random walk in one dimension is \emph{recurrent},
in that it visits the origin (and indeed \emph{any} site on $\ZZ$)
infinitely often. Moreover, it can be shown to be
\emph{null-recurrent}, which means that the expected time to
return to the origin is infinite. That is to say, return to the
origin is guaranteed, but it takes very long until this happens.

Fluctuations of the random walk can be characterized further via the
central limit theorem (CLT), which amounts to saying that the
distribution of $X_n$ is asymptotically normal, with mean $n(p-q)$
and variance $4npq$:
\begin{align}
  \lim_{n\to\infty}&\PW\left\{\frac{X_n-n(p-q)}{\sqrt{4npq}}\le x\right\}
  \notag\\
  &=
  \varPhi(x):=\frac{1}{\sqrt{2\pi}}\int_{-\infty}^x \e^{-y^2/2}\,\dif y.
  \label{Phi}
\end{align}

These results can be extended to more general walks in one
dimension, and also to higher dimensions. For instance, the
criterion of recurrence for a general one-dimensional random walk is
that it is unbiased, $\EW (X_1-X_0)=0$. In the two-dimensional case,
in addition one needs $\EW |X_1-X_0|^2<\infty$. In higher
dimensions, any random walk (which does not reduce to lower
dimension) is transient.

\subsection{Random Environments and Random\\ Walks}\label{sec:RWRE}

The definition of an RWRE involves two ingredients: (i) the
\emph{environment}, which is randomly chosen but remains fixed
throughout the time evolution, and (ii) the \emph{random walk},
whose transition probabilities are determined by the environment.
The set of environments (sample space) is denoted by
$\Omega=\{\omega\}$, and we use $\PP$ to denote the probability
distribution on this space. For each $\omega\in\Omega$, we define
the random walk in the environment $\omega$ as the
(time-homogeneous) Markov chain $\{X_t,\,t=0,1,2,\dots\}$ on
$\ZZ^d$ with certain (random) transition probabilities
\begin{equation}\label{p_xy}
  p(x,y,\omega)=\PW^\omega\{X_1=y\,|\,X_0=x\}.
\end{equation}
The probability measure $\PW^\omega$ that determines the
distribution of the random walk in a given environment $\omega$ is
referred to as the \emph{quenched} law. We often use a subindex to
indicate the initial position of the walk, so that e.g.\
$\PW^\omega_x\{X_0=x\}=1$.

By averaging the quenched probability $\PW^\omega_x$ further, with
respect to the environment distribution, we obtain the
\emph{annealed} measure $\Prob_x=\PP\times \PW^\omega_x$, which
determines the probability law of the RWRE:
\begin{equation}\label{annealed}
  \Prob_x(A)=\int_\Omega \PW_x^\omega(A)\PP(\dif\omega)
  =\EE \PW_x^\omega(A).
\end{equation}
Expectation with respect to the annealed measure $\Prob_x$ will be
denoted by $\Exp_x$.

Equation (\ref{annealed}) implies that if some property $A$ of the
RWRE holds almost surely (a.s.) with respect to the quenched law
$\PW_x^\omega$ for almost all environments (i.e., for all
$\omega\in\Omega'$ such that $\PP(\Omega')=1$), then this property
is also true with probability $1$ under the annealed law
$\Prob_x$.

Note that the random walk $X_n$ is a Markov chain only
conditionally on the fixed environment (i.e., with respect to
$\PW_x^\omega$), but the Markov property fails under the annealed
measure $\Prob_x$. This is because the past history cannot be
neglected, as it tells what information about the medium must be
taken into account when averaging with respect to environment.
That is to say, the walk learns more about the environment by
taking more steps. (This idea motivates the method of
``environment viewed from the particle'', see Section
\ref{sec:Lagrange} below.)

The simplest model is the nearest-neighbor one-dimen\-sio\-nal walk, with
transition probabilities
$$
  p(x,y,\omega)=\left\{\begin{array}{ll}
  p_x&\text{if}\ \ y=x+1,\\
  q_x&\text{if}\ \ y=x-1,\\
  0&\text{otherwise},
  \end{array}\right.
$$
where $p_x$ and $q_x=1-p_x$ ($x\in\ZZ$) are random variables on
the probability space $(\Omega,\PP)$. That is to say, given the
environment $\omega\in\Omega$, the random walk currently at point
$x\in\ZZ$ will make a one-unit step to the right, with probability
$p_x$, or to the left, with probability $q_x$. Here the
environment is determined by the sequence of random variables
$\{p_x\}$. For the most of the article, we assume that the random
probabilities $\{p_x,\,x\in\ZZ\}$ are independent and identically
distributed (i.i.d.), which is referred to as \emph{i.i.d.\
environment}. Some extensions to more general environments will be
mentioned briefly in Section \ref{sec:variants}. The study of RWRE
is simplified under the following natural condition called
\emph{(uniform) ellipticity}:
\begin{equation}\label{ellipticity}
 0<\delta\le p_x\le 1-\delta<1,\ \ x\in\ZZ,\quad \PP\text{-a.s.}
\end{equation}
which will be frequently assumed in the sequel.

\section{Transience and Recurrence}\label{sec:TvR}
In this section, we discuss a criterion for the RWRE to be transient or
recurrent.  The
following theorem is due to Solomon (1975).
\begin{theorem}\label{Solomon1}
Set $\rho_x:=q_x/p_x$, $x\in\ZZ$, and $\eta:=\EE\ln \rho_0$.

{\rm (i)} If\/ $\eta\ne 0$ then $X_t$ is transient
($\Prob_0$-a.s.); moreover, if\/ $\eta<0$ then
$\lim_{t\to\infty}X_t=+\infty$, while if\/ $\eta>0$ then\/
$\lim_{t\to\infty}X_t=-\infty$ ($\Prob_0$-a.s.).

{\rm (ii)} If\/ $\eta=0$ then $X_t$ is recurrent ($\Prob_0$-a.s.);
moreover,
$$
\varlimsup_{t\to\infty} X_t=+\infty,\quad \varliminf_{t\to\infty}
X_t=-\infty,\quad \Prob_0\text{-a.s.}
$$
\end{theorem}

Let us sketch the proof.
Consider the hitting times $T_x:=\min\{t\ge 0: X_t=x\}$ and denote by $f_{xy}$
the quenched first-passage probability from $x$ to $y$:
$$
  f_{xy}:=\PW^\omega_x\{1\le T_y<\infty\}.
$$
Starting from $0$ the first step of the walk may be either to the
right or to the left, hence by the Markov property the return
probability $f_{00}$ can be decomposed as
\begin{equation}\label{f00}
  f_{00}=p_0f_{10}+q_0f_{-1,0}.
\end{equation}
To evaluate $f_{10}$, for $n\ge 1$ set
$$
u_i\equiv u_i^{(x)}:=\PW^\omega_i\{T_0<T_x\},\quad 0\le i\le x,
$$
which is the probability to reach $0$ prior to $x$, starting from
$i$. Clearly,
\begin{equation}\label{f_k0}
  f_{10}=\lim_{x\to\infty} u_1^{(x)}. 
\end{equation}
Decomposition with respect to the first step yields the difference equation
\begin{equation}\label{u}
  u_i=p_i u_{i+1}+q_i u_{i-1}, \quad 0<i<n,
\end{equation}
with the boundary conditions
\begin{equation}\label{u0}
  u_0=1,\quad u_x=0.
\end{equation}
Using $p_x+q_x=1$, eqn (\ref{u}) can be rewritten as
$$
u_{x+1}-u_{x}=\rho_x (u_{x}-u_{x-1}),
$$
whence by iterations
\begin{equation}\label{eq:u-u}
 u_{x+1}-u_{x}=(u_1-u_0)\prod_{j=1}^x \rho_j.
\end{equation}
Summing over $x$ and using the boundary conditions (\ref{u0}) we
obtain
\begin{equation}\label{u1}
1-u_1=\Biggl(\sum_{x=0}^{n-1}\prod_{j=1}^x\rho_j\Biggr)^{-1}
\end{equation}
(if $x=0$, the product over $j$ is interpreted as $1$).
In view of eqn (\ref{f_k0}) it follows that $f_{10}=1$ if and only
if the right-hand side of eqn (\ref{u1}) tends to $0$, that is,
\begin{equation}\label{series}
  \sum_{x=1}^{\infty} \exp (Y_x)=\infty,\quad
  Y_x:=\sum_{j=1}^x \ln\rho_j.
\end{equation}
Note that the random variables $\ln\rho_j $ are i.i.d., hence
by the strong LLN
$$
  \lim_{x\to\infty}\frac{Y_x}{x}
  =\EE\ln \rho_0\equiv\eta, \quad \PP\text{-a.s.}
$$
That is, the general term of the series (\ref{series}) for
large $x$ behaves like $\exp(x\eta)$, hence for $\eta>0$ the
condition (\ref{series}) holds true (and so $f_{10}=1$), whereas
for $\eta<0$ it fails (and so $f_{10}<1$).

By interchanging the roles of $p_x$ and $q_x$, we also have $f_{-1,0}<1$ if
$\eta>0$ and $f_{-1,0}=1$ if $\eta<0$. From
eqn (\ref{f00}) it then follows that in both cases $f_{00}<1$,
i.e.\ the random walk is transient.

In the critical case, $\eta=0$, by a general result from
probability theory,
$Y_x\ge0$ for infinitely many $x$ ($\PP$-a.s.), and so the series in
eqn (\ref{series}) diverges. Hence, $f_{10}=1$ and, similarly,
$f_{-1,0}=1$, so by eqn (\ref{f00}) $f_{00}=1$, i.e.\ the random
walk is recurrent.

It may be surprising that the critical parameter appears in the
form $\eta=\EE\ln\rho_0$, as it is probably more natural to
expect, by analogy with the ordinary random walk, that the RWRE
criterion would be based on the mean drift, $\EE(p_0-q_0)$. In the
next section we will see that the sign of $d$ may be misleading.

A canonical model of RWRE is specified by the assumption that the random
variables $p_x$ take only two values, $\beta$ and $1-\beta$, with probabilities
\begin{equation}\label{example}
\PP\{p_x=\beta\}=\alpha,\quad \PP\{p_x=1-\beta\}=1-\alpha,
\end{equation}
where $0<\alpha<1$, \,$0<\beta<1$. Here $\eta=(2\alpha-1)\ln
(1+(1-2\beta)/\beta)$, and it is easy to see that, e.g., $\eta<0$
if $\alpha<\frac12$, $\beta<\frac12$ or $\alpha>\frac12$,
$\beta>\frac12$. The recurrent region where $\eta=0$ splits into
two lines, $\beta=\frac12$ and $\alpha=\frac12$. Note that the
first case is degenerate and amounts to the ordinary symmetric
random walk, while the second one (except where $\beta=\frac12$)
corresponds to Sinai's problem (see Section \ref{sec:Sinai}). A
``phase diagram'' for this model, showing various limiting regimes
as a function of the parameters $\alpha$, $\beta$, is presented in
Figure 1.

\begin{figure}[htbp]
\linethickness{0.4mm}
\hspace{-5mm}\includegraphics[width=3.0in]{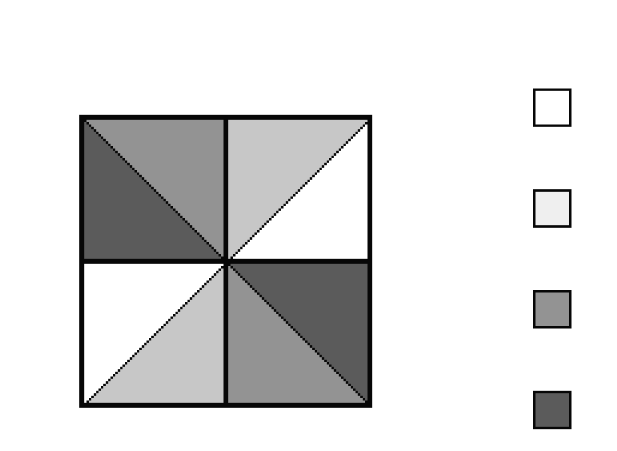}
\put(-10,124.5){\mbox{$\eta<0$}} \put(-10,112.5){\mbox{$v>0$}}
\put(-10,89){\mbox{$\eta<0$}} \put(-10,77){\mbox{$v=0$}}
\put(-10,53.5){\mbox{$\eta>0$}} \put(-10,41.5){\mbox{$v=0$}}
\put(-10,17){\mbox{$\eta>0$}} \put(-10,5){\mbox{$v<0$}}
\put(-196,7){\mbox{$0$}} \put(-197,115){\mbox{$1$}}
\put(-198,63){\mbox{$\beta$}} \put(-84,7){\mbox{$1$}}
\put(-141,6){\mbox{$\alpha$}} \put(-69,63){\mbox{$\beta=\frac12$}}
\put(-138,133){{\vector(0,-1){12}}}
\put(-151.5,137){\mbox{$\alpha=\frac12$}}
\put(-72,66.3){{\vector(-1,0){12}}} \caption{\small Phase diagram
for the canonical model, eqn (\ref{example}). In the regions where
$\eta<0$ or $\eta>0$, the RWRE is transient to $+\infty$ or
$-\infty$, respectively. The recurrent case, $\eta=0$, arises when
$\alpha=\frac12$ or $\beta=\frac12$. The asymptotic velocity
$v:=\lim_{t\to\infty} X_t/t$ is given by eqn (\ref{drift}). Adapted
from Hughes B.D. (1996) \emph{Random Walks and Random Environments.
Volume 2: Random Environments}, Ch.~6, p.~391. Clarendon Press,
Oxford, by permission of Oxford University Press.}
\end{figure}
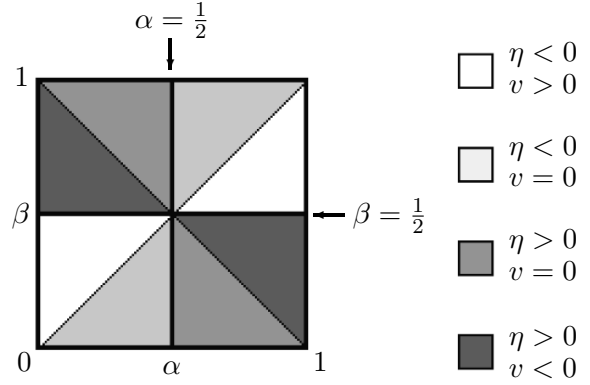

\section{Asymptotic Velocity}\label{sec:Stokes}
In the transient case the walk escapes to infinity, and it is
reasonable to ask at what speed. For a non-random environment,
$p_x\equiv p$, the answer is given by the LLN, eqn (\ref{LLN}).
For the simple RWRE, the asymptotic velocity was obtained by
Solomon (1975). Note that by Jensen's inequality,
$(\EE\rho_0)^{-1}\le \EE \rho_0^{-1}$.

\begin{theorem}\label{Solomon2}
\,The limit \,$v:=\lim_{t\to\infty} X_t/t$ \,exists\break ($\Prob_0$-a.s.) and is
given by
\begin{equation}\label{drift}
 v=\left\{
  \begin{array}{cl}
  \displaystyle
  \frac{1-\EE\rho_0}{1+\EE\rho_0}&\text{if}\
  \EE\rho_0<1,\\[.7pc]
  \displaystyle-\frac{1-\EE \rho_0^{-1}}{1+\EE \rho_0^{-1}}&\text{if}\
  \EE \rho_0^{-1}<1,\\[.7pc]
  0&\text{otherwise}.
\end{array}
\right.
\end{equation}
\end{theorem}

Thus, the RWRE has a well-defined non-zero\break asymptotic velocity
except when $(\EE\rho_0)^{-1}\le 1\le \EE\rho_0^{-1}$. For
instance, in the canonical example eqn (\ref{example}) (see Figure
1) the criterion $\EE\rho_0<1$ for the velocity $v$ to be positive
amounts to the condition that both $(1-\alpha)/\alpha$ and
$(1-\beta)/\beta$ lie on the same side of point $1$.

The key idea of the proof is to analyze the hitting times $T_n$
first, deducing results for the walk $X_t$ later. More
specifically, set $\tau_{i}=T_{i}-T_{i-1}$, which is the time to
hit $i$ after hitting $i-1$ (providing that $i>X_0$). If $X_0=0$
and $n\ge 1$ then $T_n=\tau_1+\dots+\tau_n$. Note that in fixed
environment $\omega$ the random variables $\{\tau_i\}$ are
independent, since the quenched random walk ``forgets'' its past.
Although there is no independence with respect to the annealed
probability measure $\Prob_0$, one can show that, due to the
i.i.d.\ property of the environment, the sequence $\{\tau_i\}$ is
ergodic and therefore satisfies the LLN:
$$
\frac{T_n}{n}=\frac{\tau_1+\dots+\tau_n}{n}\to\Exp_0\tau_1, \quad
\Prob_0\text{-a.s.}
$$
In turn, this implies
\begin{equation}\label{X/t}
\frac{X_t}{t}\to \frac{1}{\Exp_0 \tau_1},\quad \Prob_0\text{-a.s.}
\end{equation}
(the clue is to note that $X_{T_n}=n$).

To compute the mean value $\Exp_0\tau_1$,
observe that
\begin{equation}\label{tau_1}
\tau_{1}={\mathbf 1}_{\{X_1=1\}}+{\mathbf 1}_{\{X_1=-1\}}
(1+\tau'_{0}+\tau'_{1}),
\end{equation}
where ${\mathbf 1}_{A}$ is the indicator of event $A$ and $\tau'_0$,
$\tau'_1$ are, respectively, the times to get from $-1$ to $0$ and then from
$0$ to $1$.
Taking expectations in a fixed environment $\omega$,
we obtain
\begin{equation}\label{E-tau-0}
\EW_0^\omega\tau_{1}=p_0+q_0 (1+\EW_0^\omega\tau'_{0}+\EW_0^\omega\tau_{1}),
\end{equation}
and so
\begin{equation}\label{E-tau}
 \EW_0^\omega \tau_{1}=1+\rho_0+\rho_0 \EW_0^\omega\tau'_{0}.
\end{equation}
Note that $\EW_0^\omega\tau'_{0}$ is
a function of $\{p_x,\, x<0\}$ and hence
is independent of $\rho_0=q_0/p_0$. Averaging eqn (\ref{E-tau}) over
the environment and using $\Exp_0 \tau'_{0}=\Exp_0\tau_{1}$
yields
\begin{equation}\label{tau1}
\Exp_0\tau_{1}=\left\{
  \begin{array}{cl}
  \displaystyle
  \frac{1+\EE\rho_0}{1-\EE\rho_0}&\text{if}\
  \EE\rho_0<1,\\
  \infty&\text{if}\ \EE\rho_0\ge 1,
\end{array}
\right.
\end{equation}
and by eqn (\ref{X/t}) ``half'' of eqn (\ref{drift}) follows. The
other half, in terms of $\EE\rho_0^{-1}$, can be obtained by interchanging
the roles of $p_x$ and $q_x$, whereby $\rho_0$ is replaced with
$\rho_0^{-1}$.

Let us make a few remarks concerning Theorems \ref{Solomon1} and
\ref{Solomon2}. First of all, note that by Jensen's inequality
$\EE \ln\rho_0\le\ln\EE \rho_0$, with a strict inequality whenever
$\rho_0$ is non-degenerate. Therefore, it may be possible that,
with $\Prob_0$-probability $1$, $X_t\to\infty$ but $X_t/t\to 0$
(see Figure 1). This is quite unusual as compared to the ordinary
random walk (see Section \ref{sec:1.1}), and indicates some kind
of slowdown in the transient case.

Furthermore, by Jensen's inequality
$$
\EE\rho_0=\EE p_0^{-1} -1\ge (\EE p_0)^{-1}-1,
$$
so eqn (\ref{drift}) implies that if $\EE\rho_0<1$ then
$$
0<v\le 2\EE p_0-1=\EE(p_0-q_0),
$$
and the inequality is strict if $p_0$ is genuinely random (i.e.,
does not reduce to a constant). Hence, the asymptotic velocity $v$
is less than the mean drift $\EE (p_0-q_0)$, which is yet another
evidence of slowdown. What is even more surprising is that it is
possible to have $\EE (p_0-q_0)>0$ but $\eta=\EE\ln \rho_0>0$, so
that $\PP_0$-a.s.\ $X_t\to-\infty$ (although with velocity $v=0$).
Indeed, following Sznitman (2004) suppose that
$$
\PP\{p_0=\beta\}=\alpha,\quad \PP\{p_0=\gamma\}=1-\alpha,
$$
with $\alpha>\frac{1}{2}$. Then $\EE p_0\ge\alpha\beta>\frac12$ if
$1>\beta>\frac{1}{2\alpha}$, hence $\EE (p_0-q_0)=2\EE p_0-1>0$.
On the other hand,
$$
\EE \ln \rho_0=\alpha\ln
\frac{1-\beta}{\beta}+(1-\alpha)\ln\frac{1-\gamma}{\gamma}>0,
$$
if $\gamma$ is sufficiently small.

\section{Critical Exponent, Excursions and Traps}\label{sec:Excursions}
Extending the previous analysis of the hitting times, one can
obtain useful information about the limit distribution of $T_n$ (and hence
$X_t$). To appreciate this, note that from the recursion (\ref{tau_1}) it
follows
$$
\tau_{1}^s={\mathbf 1}_{\{X_1=1\}}+{\mathbf 1}_{\{X_1=-1\}}
(1+\tau'_{0}+\tau'_{1})^s,
$$
and, similarly to eqn (\ref{E-tau-0}),
$$
   \EW_0^\omega \tau_{1}^s=p_0+q_0 \EW_0^\omega (1+\tau'_{0}+\tau'_{1})^s.
$$
Taking here expectation $\EE$, one can deduce that
$\Exp_0\tau_1^s\allowbreak<\infty$ if and only if $\EE
\rho_0^s<1$. Therefore, it is natural to expect that the root
$\kappa$ of the equation
\begin{equation}\label{kappa}
  \EE \rho_0^\kappa=1
\end{equation}
plays the role of a critical exponent responsible for the growth
rate (and hence, for the type of the limit distribution) of the
sum $T_n=\tau_1+\dots+\tau_n$. In particular, by analogy with sums
of i.i.d.\ random variables
one can expect that if $\kappa>2$ then $T_n$ is asymptotically
normal, with the standard scaling $\sqrt{n}$, while for $\kappa<2$
the limit law of $T_n$ is stable (with index $\kappa$) under
scaling ${}\approx n^{1/\kappa}$.

Alternatively, eqn (\ref{kappa}) can be obtained from consideration of
excursions of the random walk. Let $T^{L}_{11}$ be the left excursion time from
site $1$, that is the time to return to $1$ after moving to the left at the
first step. 
If $\eta=\EE\ln \rho_0<0$, then $T^{L}_{11}<\infty$ ($\Prob_0$-a.s.). Fixing
an environment $\omega$, let $w_1=\EW_1^\omega T^{L}_{11}$ be the quenched
mean duration of the excursion $T^{L}_{11}$ and observe that
$w_1=1+\EW_0^\omega\tau_1$, where $\tau_1$ is the time to get back to $1$
after stepping to $0$.

As a matter of fact, this representation and eqn (\ref{tau1}) imply that
the annealed mean duration of the left excursion, $\Exp_0 T^{L}_{11}$,
is given by
\begin{equation}\label{T'}
\EE w_1=
\left\{\begin{array}{cl}
\displaystyle
\frac{2}{1-\EE\rho_0}&\text{if}\ \EE\rho_0<1,\\
\infty&\text{if}\ \EE\rho_0\ge 1.
\end{array}\right.
\end{equation}
Note that in the latter case (and bearing in mind $\eta<0$), the random walk
starting from $1$ will eventually drift to $+\infty$, thus making only a finite
number of visits to $0$, but the expected number of such visits is infinite.

In fact, our goal here is to characterize the distribution of
$w_1$ under the law $\PP$. To this end, observe that the excursion
$T^{L}_{11}$ involves at least two steps (the first and the last
ones) and, possibly, several left excursions from $0$, each with
mean time $w_0=\EW_0^\omega T^{L}_{00}$. Therefore,
\begin{equation}\label{w1}
 w_1=2+\sum_{j=1}^\infty q_0^jp_0 (j w_0) = 2+\rho_0 w_0.
\end{equation}
By the translation invariance of the environment,
the random variables $w_1$ and $w_0$ have the same distribution.
Furthermore, similarly to recursion (\ref{w1}), we have
$w_0=2+\rho_{-1} w_{-1}$. This implies that $w_0$ is a function of
$p_x$ with $x\le -1$ only, and hence $w_0$ and $\rho_0$ are
independent random variables. Introducing the Laplace transform
$\phi(s)=\EE \exp(-s w_1)$ and conditioning on $\rho_0$, from eqn
(\ref{w1}) we get the equation
\begin{equation}\label{phi}
\phi(s)=\e^{-2s}\EE \phi(s\rho_0).
\end{equation}
Suppose that
$$
  1-\phi(s)\sim as^\kappa,\quad s\to 0,
$$
then eqn (\ref{phi}) amounts to
$$
1-as^\kappa+\cdots=(1-2s+\cdots)(1-as^\kappa
\EE\rho_0^\kappa+\cdots).
$$
Expanding the product on the right, one can see that a solution with $\kappa=1$
is possible only if $\EE\rho_0<1$, in which case
$$
  a=\EE w_1 = \frac{2}{1-\EE\rho_0}.
$$
We have already obtained this result in eqn (\ref{T'}).

The case $\kappa<1$ is possible if $\EE\rho_0^\kappa=1$, which is
exactly eqn (\ref{kappa}). Returning to $w_1$, one expects a slow
decay of the distribution tail,
$$
  \PP\{w_1>t\}\sim b\, t^{-1/\kappa},\quad t\to\infty.
$$
In particular, in this case the annealed mean duration of the left excursion
appears to be infinite.

Although the above considerations point to the critical parameter
$\kappa$, eqn (\ref{kappa}), which may be expected to determine
the slowdown scale, they provide little explanation of a mechanism
of the slowdown phenomenon. Heuristically, it is natural to
attribute the slowdown effects to the presence of \emph{traps} in
the environment, which may be thought of as regions that are easy
to enter but hard to leave. In the one-dimensional case, such a
trap would occur, for example, between two long series of
successive sites where the probabilities $p_x$ are fairly large
(on the left) and small (on the right).

Remarkably, traps can be characterized quantitatively with regard
to the properties of the random environment, by linking them to
certain large deviation effects (see Sznitman (2002, 2004)). The
key role in this analysis is played by the function
$F(u):=\ln\EE\rho_0^u$, \,$u\in\RR$. Suppose that
$\eta=\EE\ln\rho_0<0$ (so that by Theorem \ref{Solomon1} the RWRE
tends to $+\infty$, $\Prob_0$-a.s.) and also that $\EE\rho_0>1$
and $\EE\rho_0^{-1}>1$ (so that by  Theorem \ref{Solomon2},
$v=0$). The latter means that $F(1)>0$ and $F(-1)>0$, and since
$F$ is a smooth strictly convex function and $F(0)=0$, it follows
that there is the second root $0<\kappa<1$, so that $F(\kappa)=0$,
i.e., $\EE\rho_0^\kappa=1$ (cf.\ eqn (\ref{kappa})).

Let us estimate the probability to have a trap in $U=[-L,L]$ where the
RWRE will spend anomalously long time. Using eqn (\ref{u1}), observe that
$$
  \PW_1^\omega \{T_0<T_{L+1}\}\ge 1-\exp\{-L S_L\},
$$
where $S_L:=L^{-1}\sum_{x=1}^L\ln\rho_x\to\eta<0$ as $L\to\infty$. However,
due to large deviations $S_L$ may exceed level $\epsilon>0$ with probability
$$
  \PP\{S_L>\epsilon\}\sim\exp\{-L I(\epsilon)\},\quad L\to\infty,
$$
where $I(x):=\sup_u\{ux-F(u)\}$ is the Legendre transform of $F$. We
can optimize this estimate by assuming that $\epsilon L\ge \ln n$
and minimizing the ratio $I(\epsilon)/\epsilon$. Note that $F(u)$
can be expressed via the inverse Legendre transform,
$F(u)=\sup_x\{xu-I(x)\}$, and it is easy to see that if
$\kappa:=\min_{\epsilon>0} I(\epsilon)/\epsilon$ then $F(\kappa)=0$,
so $\kappa$ is the second (positive) root of~$F$.

The ``left'' probability $\PW_{-1}^\omega \{T_0<T_{-L-1}\}$ is estimated
in a similar fashion, and one can deduce that for some constants $K>0$, $c>0$
and any $\kappa'>\kappa$, for large $n$
$$
  \PP\Bigl\{\PW_0^\omega \Bigl\{\max_{k\le n} |X_k|\le K\ln n\Bigr\}\ge c\Bigr\}\ge
n^{-\kappa'}.
$$
That is to say, this is a bound on the probability to see a trap
centered at $0$, of size $\approx \ln n$, which will retain the
RWRE for at least time $n$. It can be shown that, typically, there
will be many such traps both in $[-n^{\kappa'},0]$ and
$[0,n^{\kappa'}]$, which will essentially prevent the RWRE from
moving at distance $n^{\kappa'}$ from the origin before time $n$.
In particular, it follows that $\lim_{n\to\infty}
X_n/n^{\kappa'}=0$ for any $\kappa'>\kappa$, so recalling that
$0<\kappa<1$, we have indeed a sublinear growth of $X_n$. This
result is more informative as compared to Theorem \ref{Solomon2}
(the case $v=0$), and it clarifies the role of traps (see more
details in Sznitman (2004)). The non-trivial behavior of the RWRE
on the precise growth scale, $n^\kappa$, is characterized in the
next section.

\section{Limit Distributions}\label{sec:KKS}
Considerations in Section \ref{sec:Excursions} suggest that the
exponent $\kappa$, defined as the solution of eqn (\ref{kappa}),
characterizes environments in terms of duration of left
excursions. These heuristic arguments are confirmed by a limit
theorem by Kesten \emph{et al.}\ (1975), which specifies the
slowdown scale. We state here the most striking part of their
result. Denote $\ln^+\!u:=\max\{\ln u,0\}$; by an \emph{arithmetic
distribution} one means a probability law on $\RR$ concentrated on
the set of points of the form $0$, $\pm c$, $\pm2 c$, $\dots$

\begin{theorem}\label{KKS} Assume that\/ $-\infty\le\eta=\EE\ln\rho_0<0$
and the distribution of\/ $\ln\rho_0$ is non-arithmetic (excluding a possible
atom at $-\infty$). Suppose that the root $\kappa$ of equation (\ref{kappa})
is such that $0<\kappa<1$ and $\EE \rho_0^\kappa\ln^+\!\rho_0<\infty$.
Then
\begin{gather*}
\lim_{n\to\infty}\Prob_0\{n^{-1/\kappa}\,T_n\le t\}=L_\kappa
(t),\\
\lim_{t\to\infty}\Prob_0\{t^{-\kappa}X_t\le x\}=1-L_\kappa
(x^{-1/\kappa}),
\end{gather*}
where $L_\kappa(\cdot)$ is the distribution function of a stable
law with index $\kappa$, concentrated on $[0,\infty)$.
\end{theorem}

General information on stable laws can be found in many probability
books;
we only mention here that the Laplace transform of a stable
distribution on $[0,\infty)$ with index $\kappa$ has the form
$\phi(s)=\exp\{-Cs^\kappa\}$.

Kesten \emph{et al.}\ (1975) also consider the case $\kappa\ge 1$.
Note that for $\kappa>1$, we have
$\EE\rho_0<(\EE\rho_0^\kappa)^{1/\kappa}=1$, so $v>0$ by eqn (14).
For example, if $\kappa>2$ then, as expected (see
Section~\ref{sec:Excursions}),
\begin{gather*}
\lim_{n\to\infty}\Prob_0\biggl\{\frac{T_n-n/v}{\sigma\sqrt{n}}
\le t\biggr\}=\varPhi(t),\\[.5pc]
\lim_{t\to\infty}\Prob_0\biggl\{\frac{X_t-tv}{v^{3/2}\sigma\sqrt{t}}
\le x\biggr\}=\varPhi(x).
\end{gather*}

Let us describe an elegant idea of the proof based on a suitable
renewal structure. (i) Let $U_i^n$ ($i\le n$) be the number of
left excursions starting from $i$ up to time $T_n$, and note that
$T_n=n+2\sum_{i} U_i^n$. Since the walk is transient to $+\infty$,
the sum $\sum_{i\le 0} U_i^n$ is finite ($\Prob_0$-a.s.) and so
does not affect the limit. (ii) Observe that if the environment
$\omega$ is fixed then the conditional distribution of $U_j^n$,
given $U_{j+1}^n,\dots,U_n^n=0$, is the same as the distribution
of the sum of $1+U_{j+1}^n$ i.i.d.\ random variables
$V_1,V_2,\dots$, each with geometric distribution
$\PW_0^\omega\{V_i=k\}=p_j q_j^k$ ($k=0,1,2,\dots$). Therefore,
the sum $\sum_{i=1}^n U_i^n$ (read from right to left) can be
represented as $\sum_{t=0}^{n-1} Z_t$, where $Z_0=0,Z_1,
Z_2,\dots$ is a branching process (in random environment
$\{p_j\}$) with one immigrant at each step and the geometric
offspring distribution with parameter $p_j$ for each particle
present at time $j$. (iii) Consider the successive
``regeneration'' times $\tau_k^*$, at which the process $Z_t$
vanishes. The partial sums $W_k:=\sum_{\tau_k^*\le
t<\tau^*_{k+1}}Z_t$ form an i.i.d.\ sequence, and the proof
amounts to showing that the sum of $W_k$ has a stable limit of
index $\kappa$. (iv) Finally, the distribution of $W_0$ can be
approximated using $M_0:=\sum_{t=1}^\infty\prod_{j=0}^{n-1}\rho_j$
(cf.\ eqn (\ref{u1})), which is the quenched mean number of total
progeny of the immigrant at time $t=0$. Using Kesten's renewal
theorem, it can be checked that $\PP\{M_0>x\}\sim K x^{-\kappa}$
as $x\to\infty$, so $M_0$ is in the domain of attraction of a
stable law with index $\kappa$, and the result follows.

Let us emphasize the significance of the regeneration times
$\tau_i^*$. Returning to the original random walk, one can see
that these are times at which the RWRE hits a new ``record'' on
its way to $+\infty$, never to backtrack again. The same idea
plays a crucial role in the analysis of the RWRE in higher
dimensions (see Sections \ref{sec:0-1}, \ref{sec:Kalikow} below).

Finally, note that the condition $-\infty\le\eta<0$ allows
$\PP\{p_0=1\}>0$, so the distribution of $\rho_0$ may have an atom
at $0$ (and hence $\ln\rho_0$ at $-\infty$). In view of eqn
(\ref{kappa}), no atom is possible at $+\infty$. The restriction
for the distribution of $\ln\rho_0$ to be non-arithmetic is
important. This will be illustrated in Section \ref{sec:Diodes}
where we discuss the model of random diodes.

\section{Sinai's Localization}\label{sec:Sinai}
The results discussed in Section \ref{sec:KKS} indicate that the
less transient the RWRE is (i.e., the critical exponent decreasing
to zero), the slower it moves. Sinai (1982) proved a remarkable
theorem showing that for the \emph{recurrent} RWRE (i.e., with
$\eta=\EE\ln\rho_0=0$), the slowdown effect is exhibited in a
striking way.
\begin{theorem}\label{th:Sinai}
Suppose that the environment $\{p_x\}$ is i.i.d.\ and elliptic, eqn
(\ref{ellipticity}), and assume that $\EE\ln\rho_0=0$, with
$\PP\{\rho_0=1\}<1$. Denote\/ $\sigma^2\!:=\EE \ln^2\!\rho_0$,
$0<\sigma^2<\infty$. Then there exists a function $W_n=W_n(\omega)$
of the random environment such that for any $\varepsilon>0$
\begin{equation}\label{eq:Sinai}
\lim_{n\to\infty}\Prob_0\left\{\left|\frac{\sigma^2 X_n}{\ln^2
n}-W_n\right|>\varepsilon\right\} = 0. 
\end{equation}
Moreover, $W_n$ has a limit distribution:
\begin{equation}\label{G}
\lim_{n\to\infty}\PP\left\{W_n\le x\right\}=G(x),
\end{equation}
and thus also the distribution of $\sigma^2X_n/\ln^2 n$ under
$\Prob_0$ converges to the same distribution $G(x)$.
\end{theorem}

Sinai's theorem shows that in the recurrent case, the RWRE
considered on the spatial scale $\ln^2 n$ becomes localized near
some random point (depending on the environment only). This
phenomenon, frequently referred to as \emph{Sinai's localization},
indicates an extremely strong slowdown of the motion as compared
with the ordinary diffusive behavior.

Following R\'ev\'esz (1990), let us explain heuristically why $X_n$ is measured
on the scale $\ln^2 n$. Rewrite eqn (\ref{u1}) as
\begin{equation}\label{Y1}
\PW^\omega_1
\{T_n<T_0\}
=\biggl(1+\sum_{x=1}^{n-1} 
\exp (Y_x)\biggr)^{-1},
\end{equation}
where $Y_x$ is defined in (\ref{series}). By the central limit
theorem, the typical size of $|Y_x|$ for large $x$ is of order of
$\sqrt{x}$, and so eqn (\ref{Y1}) yields
$$
 \PW_1^\omega\{T_n<T_0\}\approx \exp\{-\sqrt{n}\,\}.
$$
This suggests that the walk started at site $1$ will make about
$\exp\{\sqrt{n}\,\}$ visits to the origin before reaching level $n$.
Therefore, the first passage to site $n$ takes at least time
$\approx\exp\{\sqrt{n}\,\}$. In other words, one may expect that a
typical displacement after $n$ steps will be of order of $\ln^2 n$
(cf.\ eqn (\ref{eq:Sinai})). This argument also indicates, in the
spirit of the trapping mechanism of slowdown discussed at the end of
Section \ref{sec:Excursions}, that there is typically a trap of size
$\approx\ln^2 n$, which retains the RWRE until time $n$.

It has been shown (independently by H.~Kesten and A.O.~Golosov)
that the limit in (\ref{G}) coincides with the distribution of a
certain functional of the standard Brownian motion, with the density
function
$$
G'(x)=\frac{2}{\pi}\sum_{k=0}^\infty
\frac{(-1)^k}{2k+1}\,\exp\left\{-\frac{(2k+1)^2\pi^2}{8}\,|x|\right\}.
$$

\section{Environment Viewed from the\\ Particle}\label{sec:Lagrange}
This important technique, dating back to Kozlov and Molchanov
(1984),
has proved to be quite efficient in the study of random motions in
random media. The basic idea is to focus on the evolution of the
environment viewed from the current position of the walk.

Let $\theta$ be the shift operator acting on the space of
environments $\Omega=\{\omega\}$ as follows:
$$
\omega=\{p_x\}\stackrel{\theta}{\mapsto} \bar\omega=\{p_{x-1}\}.
$$
Consider the process
$$
\omega_n:=\theta^{X_n}\omega,\quad \omega_0=\omega,
$$
which describes the state of the environment from the point of view
of an observer moving along with the random walk $X_n$. One can show
that $\omega_n$ is a Markov chain (with respect to both
$\PW_0^\omega$ and $\Prob_0$), with the transition kernel
\begin{equation}\label{eq:T}
T(\omega,\dif\omega')=p_0\,
\delta_{\theta\omega}(\dif\omega')+q_0\,\delta_{\theta^{-1}\omega}(\dif\omega')
\end{equation}
and the respective initial law $\delta_\omega$ or $\PP$ (here
$\delta_\omega$ is the Dirac measure, i.e., unit mass at
$\omega$).

This fact as it stands may not seem to be of any practical use,
since the state space of this Markov chain is very complex.
However, the great advantage is that one can find 
an explicit invariant probability $\QQ$ for the kernel $T$ (i.e.,
such that $\QQ T=\QQ$), which is absolutely continuous with
respect to $\PP$.

More specifically, assume that $\EE\rho_0<1$ and set
$\QQ=f(\omega) \PP$, where (cf.\ eqn (\ref{drift}))
\begin{equation}\label{f}
f=v\,(1+\rho_0)\sum_{x=0}^\infty\prod_{j=1}^x\rho_j,\quad
v=\frac{1-\EE\rho_0}{1+\EE\rho_0}.
\end{equation}
Using independence of $\{\rho_x\}$, we note
$$
\int_{\Omega} \QQ(\dif \omega)=\EE
f=(1-\EE\rho_0)\sum_{x=0}^\infty (\EE\rho_0)^x=1,
$$
hence $\QQ$ is a probability measure on $\Omega$. Furthermore, for
any bounded measurable function $g$ on $\Omega$ we have
\begin{equation}\label{QT}
\begin{aligned} \QQ Tg&=\int_\Omega
Tg(\omega)\QQ(\dif\omega)=\EE f
Tg\\
&=\EE\bigl\{f\bigl[p_0\,
(g\circ \theta)+q_0\,(g\circ \theta^{-1})\bigr]\bigr\}\\
&= \EE\bigl\{g\bigl[(p_0f)\circ \theta^{-1}+(q_0 f)\circ \theta)
\bigr].
\end{aligned}
\end{equation}
By eqn (\ref{f}),
\begin{gather*}
(p_0f)\circ \theta^{-1}=v p_{-1}(1+\rho_{-1})
\sum_{x=0}^\infty\prod_{j=1}^{x}\rho_{j-1}\\
=v \biggl(1+\rho_0\sum_{x=0}^\infty\prod_{j=1}^{x}\rho_{j}\biggr)
=v+\frac{\rho_0}{1+\rho_0}\, f,
\end{gather*}
and similarly
$$
(q_0 f)\circ \theta =-v+\frac{1}{1+\rho_0}\,f.
$$
So from eqn (\ref{QT}) we obtain
$$
\QQ Tg=\EE(gf) =\int_{\Omega}g(\omega)\QQ(\dif\omega)=\QQ g,
$$
which proves the invariance of $\QQ$.

To illustrate the environment method, let us sketch the proof
of Solomon's result on the asymptotic velocity (see Theorem
\ref{Solomon2} in Section~\ref{sec:Stokes}). Set
$d(x,\omega):=\EW_x^\omega(X_1-X_0)=p_x-q_x$. Noting that
$d(x,\omega)=d(0,\theta^x\omega)$, define
$$
D_n:=\sum_{i=1}^n d(X_{i-1},\omega)=\sum_{i=1}^n
d(0,\theta^{X_{i-1}}\omega).
$$
Due to the Markov property, the process $M_n:=X_n-D_n$ is a
martingale with respect to the natural filtration
$\calF_n=\sigma\{X_1,\dots,X_n\}$ and the law $\PW_0^\omega$,
$$
\EW_0^\omega [M_{n+1}\,|\,\calF_{n}]=M_{n}\quad (\PW_0^\omega\text{-a.s.}),
$$
and it has bounded jumps, $|M_n-M_{n-1}|\le 2$. By general
results, this implies $M_n/n\to0$ ($\PW_0^\omega$-a.s.).

On the other hand, by Birkhoff's ergodic theorem
$$
\lim_{n\to\infty}\frac{D_n}{n}=\int_{\Omega}
d(0,\omega)\QQ(\dif\omega),\quad \Prob_0\text{-a.s.}
$$
The last integral is easily evaluated to yield
\begin{gather*}
\EE
(p_0-q_0)f=v\EE\sum_{x=0}^\infty\prod_{j=1}^x\rho_j(1-\rho_0) \\
=v (1-\EE\rho_0) \sum_{x=0}^\infty
(\EE\rho_0)^x=v, 
\end{gather*}
and the first part of the formula (\ref{drift}) follows.

The case $\EE\rho_0\ge 1$ can be handled using a comparison
argument (Sznitman 2004). Observe that if $p_x\le \tilde p_x$ for
all $x$ then for the corresponding random walks we have $X_t\le
\tilde X_t$ ($\PW_0^\omega$-a.s.). We now define a suitable
dominating random medium by setting (for $\gamma>0$)
$$
\tilde p_x:=\frac{p_x}{1+\gamma}+\frac{\gamma}{1+\gamma}\ge p_x.
$$
Then $\EE\tilde \rho_0= \EE q_0/(p_0+\gamma)<1$ if $\gamma$ is
large enough, so by the first part of the theorem,
$\PW_0^\omega$-a.s.,
\begin{equation}\label{S1}
\varlimsup_{n\to\infty}\frac{X_n}{n}\le
\lim_{n\to\infty}\frac{\tilde
X_n}{n}=\frac{1-\EE\tilde\rho_0}{1+\EE\tilde\rho_0}.
\end{equation}
Note that $\EE\tilde\rho_0$ is a continuous function of $\gamma$
with values in $[0,\EE\rho_0]\ni1$, so there exists $\gamma^*$
such that $\EE\tilde\rho_0$ attains the value $1$. Passing to the
limit in eqn (\ref{S1}) as $\gamma\uparrow\gamma^*$, we obtain
$\varlimsup_{n\to\infty} X_n/n\le0$ ($\PW_0^\omega$-a.s.).
Similarly, we get the reverse inequality, which proves the second
part of the theorem.

A more prominent advantage of the environment method is that it
naturally leads to statements of CLT type. A key step
is to find a function $H(x,t,\omega)=x-vt+h(x,\omega)$
(called \emph{harmonic coordinate}) such that the process
$H(X_n,n,\omega)$ is a martingale. To this end, by the Markov
property it suffices to have
$$
  \EW_{X_n}^\omega H(X_{n+1},n+1,\omega)=H(X_n,n,\omega),\quad \PW_0^\omega\text{-a.s.}
$$
For $\Delta(x,\omega):=h(x+1,\omega)-h(x,\omega)$ this condition
leads to the equation
$$
 \Delta(x,\omega)=\rho_x\Delta(x-1,\omega)+v-1+(1+v)\rho_x.
$$
If $\EE\rho_0<1$ (so that $v>0$), there exists a bounded solution
$$
 \Delta(x,\omega)=v-1+2v\sum_{k=0}^\infty\prod_{i=0}^k \rho_{x-i},
$$
and we note that $\Delta(x,\omega)=\Delta(0,\theta^x\omega)$ is a stationary sequence with mean
$\EE\Delta(x,\omega)=0$. Finally, setting $h(0,\omega)=0$ we find
$$
  h(x,\omega)=\left\{\begin{array}{ll}
  \displaystyle
  \hphantom{-}\sum_{k=0}^{x-1}\Delta(k,\omega),&x>0,\\
\displaystyle
 - \sum_{k=1}^{-x}\Delta(-k,\omega),&x<0.
 \end{array}
 \right.
$$

As a result, we have the representation
\begin{equation}\label{eq:X-nv}
  X_n-nv=H(X_n,n,\omega)+h(X_n,\omega).
\end{equation}
For a fixed $\omega$, one can apply a suitable CLT for martingale differences
to the martingale term in (\ref{eq:X-nv}), while using that $X_n\sim nv$
($\Prob_0$-a.s.), the second term in (\ref{eq:X-nv}) is approximated
by the sum $\sum_{k=0}^{nv}\Delta(k,\omega)$, which can be handled via a CLT
for stationary sequences. This way, we arrive at the following result.
\begin{theorem}\label{CLT-env}
Suppose that the environment is elliptic, eqn (\ref{ellipticity}),
and such that $\EE\rho_0^{2+\varepsilon}<1$ for some $\varepsilon
>0$ (which implies that $\EE\rho_0<1$ and hence $v>0$). Then there
exists a non-random $\sigma^2>0$ such that
$$
\lim_{n\to\infty}\Prob_0\left\{\frac{X_n-nv}{\sqrt{n\sigma^2}}\le
x\right\}=\varPhi(x).
$$
\end{theorem}

Note that this theorem is parallel to the result by Kesten
\emph{et al.}\ (1975) on asymptotic normality when $\kappa>2$ (see
Section~\ref{sec:KKS}). The assumptions in Theorem \ref{CLT-env}
as stated are a bit more restrictive than in Theorem \ref{KKS},
but they can be relaxed. More importantly, the environment method
proves to be quite efficient in more general situations, including
non-i.i.d.\ environments and higher dimensions (at least in some
cases, e.g., for random bonds RWRE
and balanced RWRE.

\section{Diode Model}\label{sec:Diodes}
In the preceding sections (except in Section \ref{sec:KKS}, where however
we were limited to a non-arithmetic case), we assumed that
$0<p_x<1$ and therefore excluded the situation where there are
sites
through which motion is permitted in one direction only. Allowing
for such a possibility leads to the \emph{diode model} (Solomon
1975). Specifically, suppose that
\begin{equation}\label{eq:diodes}
  \PP\{p_x=\beta\}=\alpha,\quad \PP\{p_x=1\}=1-\alpha,
\end{equation}
with $0<\alpha<1$, \,$0<\beta<1$, so that with probability $\alpha$
a point $x\in\ZZ$ is a usual two-way site and with probability
$1-\alpha$ it is a repelling barrier (``diode''), through which
passage is only possible from left to right. This is an interesting
example of statistically inhomogeneous medium, where the particle
motion is strongly irreversible due to the presence of special
semi-penetrable nodes. The principal mathematical advantage of such
a model is that the random walk can be decomposed into independent
excursions from one diode to the next.

Due to diodes the random walk will eventually drift to $+\infty$.
If $\beta>\frac12$, then on average it moves faster than in a
non-random environment with $p_x\equiv \beta$. The situation where
$\beta\le\frac12$ is potentially more interesting, as then there
is a competition between the local drift of the walk to the left
(in ordinary sites) and the presence of repelling diodes on its
way. Note that $\EE\rho_0=\alpha \rho$, where
$\rho:=(1-\beta)/\beta$, so the condition $\EE\rho_0<1$ amounts to
$\beta>\alpha/(1+\alpha)$. In this case (which includes
$\beta>\frac12$), formula (\ref{drift}) for the asymptotic
velocity applies.

As explained in Section \ref{sec:Excursions}, the quenched mean
duration $w$ of the left excursion has Laplace transform given by
eqn (\ref{phi}), which now reads
$$
  \phi(s)=\e^{-2s}\bigl\{1-\alpha+\alpha\,\phi(s\rho)\bigr\}.
$$
This equation is easily solved by iterations:
\begin{equation}\label{w}
\begin{gathered}
\phi(s)=(1-\alpha)\sum_{k=0}^\infty \alpha^k \e^{-s t_k},\\
t_k:=2\sum_{j=0}^k \rho^j,
\end{gathered}
\end{equation}
hence the distribution of $w$ is given by
$$
  \PP\{w=t_k\}=(1-\alpha)\,\alpha^k,\quad k=0,1,\dots
$$
This result has a transparent probabilistic meaning. In fact, the
factor $(1-\alpha)\,\alpha^k$ is the probability that the nearest
diode on the left of the starting point occurs at distance $k+1$,
whereas $t_k$ is the corresponding mean excursion time. Note that
formula (\ref{w}) for $t_k$ easily follows from the recursion
$t_k=2+\rho t_{k-1}$ (cf.\ eqn (\ref{w1})) with the boundary
condition $t_0=2$.

A self-similar hierarchy of time scales (\ref{w}) indicates that
the process will exhibit temporal oscillations. Indeed, for
$\alpha\rho>1$ the average waiting time until passing through a
valley of ordinary sites of length $k$ is asymptotically
proportional to $t_k\sim 2\rho^k$, so one may expect the annealed
mean displacement $\Exp_0 X_t$ to have a local minimum at
$t\approx t_k$. Passing to logarithms, we note that $\ln
t_{k+1}-\ln t_k\sim \ln\rho$, which suggests the occurrence of
persistent oscillations on the logarithmic time scale, with period
$\ln\rho$. This was confirmed by Bernasconi and Schneider (1985)
who showed that for $\alpha\rho>1$
\begin{equation}\label{eq:periodic}
  \Exp_0 X_n \sim n^\kappa F(\ln n),\quad n\to\infty,
\end{equation}
where $\kappa=-\ln\alpha/\ln \rho<1$ is the solution of eqn
(\ref{kappa}) and the function $F$ is periodic with period
$\ln\rho$ (see Figure 2).

\begin{figure}[htbp]
\hspace{-2mm}\includegraphics[width=3.5in]{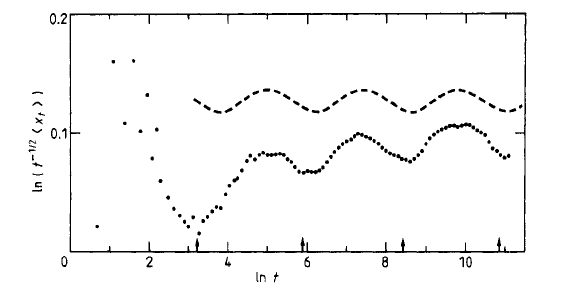}
\caption{\small Temporal oscillations for the diode model, eqn
(\ref{eq:diodes}). Here $\alpha=0.3$ and $\rho=1/0.09$, so that
$\alpha\rho>1$ and $\kappa=\frac12$. The dots represent an average
of Monte Carlo simulations over $10\,000$ samples of the environment
with a random walk of 200\,000 steps in each realization. The broken
curve refers to the exact asymptotic solution (\ref{eq:periodic}).
The arrows indicate the simulated locations of the minima $t_k$, the
asymptotic spacing of which is predicted to be $\ln \rho\approx
2.41$. Reproduced from Bernasconi J.\ and Schneider W.R.\ (1982)
Diffusion on a one-dimensional lattice with random asymmetric
transition rates. \emph{Journal of Physics A: Mathematical and
General}, {\bf 15}, L729--L734, by permission of IOP Publishing
Ltd.}
\end{figure}

In contrast, for $\alpha\rho=1$ one has
$$
  \Exp_0 X_n\sim \frac{n\ln\rho}{2\ln n},\quad n\to\infty,
$$
and there are no oscillations of the above kind.

These results illuminate the earlier analysis of the diode model
by Solomon (1975), which in the main has revealed the following.
If $\alpha\rho=1$ then $X_n$ satisfies the strong LLN:
$$
\lim_{n\to\infty} \frac{X_n}{n/\ln n}=
 \frac{\ln\rho}{2},\quad \Prob_0\text{-a.s.},
$$
while in the case $\alpha\rho>1$ the asymptotic behavior of $X_n$
is quite complicated and unusual: if $n_i\to\infty$ is a sequence
of integers such that $\{\ln n_i\}\to\gamma$ (here $\{a\}=a-[a]$
denotes the fractional part of $a$), then
the distribution of $n_i^{-\kappa} X_{n_i}$ under $\Prob_0$
converges to a non-degenerate distribution \emph{which depends on
$\gamma$.} Thus, the very existence of the limiting distribution of $X_n$ and
the limit itself heavily depend on the subsequence $n_i$ chosen to
approach infinity.

This should be compared with a more ``regular'' result in
Theorem \ref{KKS}. Note that almost all the conditions of this
theorem are satisfied in the diode model, except that here the
distribution of $\ln\rho_0$ is arithmetic (recall that the value
$\ln\rho_0=-\infty$ is permissible), so it is the discreteness of
the environment distribution that does not provide enough
``mixing'' and hence leads to such peculiar features of the
asymptotics.

\section{Some Generalizations and\\ Variations}\label{sec:variants}
Most of the results discussed above in the simplest context of
RWRE with nearest-neighbor jumps in an i.i.d.\ random environment,
have been extended to some other cases. One natural generalization
is to relax the i.i.d.\ assumption, e.g.\ by considering
stationary ergodic environments (see details in Zeitouni (2004)).
In this context, one relies on an ergodic theorem instead of the
usual strong LLN. For instance, this way one readily obtains an
extension of Solomon's criterion of transience vs.\ recurrence
(see Theorem \ref{Solomon1}, Section \ref{sec:TvR}). Other
examples include an LLN (along with a formula for the asymptotic
velocity, cf.\ Theorem \ref{Solomon2}, Section \ref{sec:Stokes}),
a CLT and stable laws for the asymptotic distribution of $X_n$
(cf.\ Theorem \ref{KKS}, Section \ref{sec:KKS}), and Sinai's
localization result for the recurrent RWRE (cf.\ Theorem
\ref{th:Sinai}, Section \ref{sec:Sinai}). Usually, however,
ergodic theorems cannot be applied directly (like, e.g., to $X_n$,
as the sequence $X_n-X_{n-1}$ is not stationary). In this case,
one rather uses the hitting times which possess the desired
stationarity (cf.\ Sections \ref{sec:Stokes},
\ref{sec:Excursions}). In some situations, in addition to
stationarity one needs suitable mixing conditions in order to
ensure enough decoupling (e.g., in Sinai's problem). The method of
environment viewed from the particle (see Section
\ref{sec:Lagrange}) is also suited very well to dealing with
stationarity.

In the remainder of this section, we describe some other
generalizations including RWRE with bounded jumps,
RWRE where randomness is attached to bonds rather than sites,
and continuous-time (symmetric) RWRE driven by the randomized
master equation.

\subsection{RWRE with Bounded Jumps}\label{sec:bounded}
The previous discussion was restricted to the case of RWRE with
nearest-neighbor jumps. A natural extension is RWRE with
\emph{bounded} jumps. Let $L,R$ be fixed natural numbers, and
suppose that from each site $x\in\ZZ$ jumps are only possible to
the sites $x+i$, \,$i=-L,\dots,R$, with (random) probabilities
\begin{gather}
\label{eq:=1}
  p_x(i)\ge 0,\quad \sum_{i=-L}^R p_x(i)=1.
\end{gather}
We assume that the random vectors $p_x(\cdot)$ determining the
environment are i.i.d.\ for different $x\in\ZZ$ (although many
results can be extended to the stationary ergodic case).

The study of asymptotic properties of such a model is essentially
more complex, as it involves products of certain random matrices
and hence must use extensively the theory of Lyapunov exponents
(see details and further references in
Br\'emont (2004)). Lyapunov exponents, being natural analogs of
logarithms of eigenvalues, characterize the asymptotic action of
the product of random matrices along (random) principal
directions, as described by Oseledec's multiplicative ergodic
theorem. In most situations, however, the Lyapunov spectrum can
only be accessed implicitly, which makes the analysis rather hard.

To explain how random matrices arise here, let us first consider a
particular case $R=1$, \,$L\ge 1$. Assume that
$p_x(-L),p_x(1)\ge\delta>0$ for all $x\in\ZZ$ (ellipticity
condition, cf.\ eqn (\ref{ellipticity})), and consider the hitting
probabilities $u_n:=\PW_n^\omega\{T_0<\infty\}$, where
$T_0:=\min\{t\ge0:X_t\le 0\}$ (cf.\ Section \ref{sec:TvR}). By
decomposing with respect to the first step, for $n\ge 1$ we obtain
the difference equation
\begin{equation}\label{eq:uLR}
 u_n=p_n(1)\,u_{n+1}+\sum_{i=0}^L p_n(-i)\, u_{n-i}
\end{equation}
with the boundary conditions $u_{0}=\dots=u_{-L+1}=1$. Using that
$1=p_n(1)+\sum_{i=0}^L p_n(-i)$, we can rewrite eqn (\ref{eq:uLR})
as
$$
p_n(1)\left(u_n-u_{n+1}\right)=\sum_{i=1}^{L}
p_n(-i)\left(u_{n-i}-u_{n}\right),
$$
or equivalently
\begin{equation}\label{v_n}
v_{n}=\sum_{i=1}^{L} b_n(i)\,v_{n-i},
\end{equation}
where $v_i:=u_{i}-u_{i+1}$ and
\begin{equation}\label{eq:b}
b_n(i):=\frac{p_n(-i)+\dots+p_n(-L)}{p_n(1)}.
\end{equation}
Recursion (\ref{v_n}) can be written in a matrix form, $V_n=M_n
V_{n-1}$, where $V_n:=(v_{n},\dots,v_{n-L+1})^\top$,
\begin{equation}\label{eq:M}
\arraycolsep=.15pc
  M_n:=\left(\begin{array}{cccc}
 b_n(1)&\dots&\dots
 &b_n(L)\\[.3pc]
 1&\dots&0&0\\[-.1pc]
 \vdots&\ddots &\vdots&\vdots\\[.1pc]
 0&\dots&1&0
 \end{array}\right),
\end{equation}
and by iterations we get (cf.\ eqn (\ref{eq:u-u}))
$$
  V_n=M_n \cdots M_1 V_0,\quad V_0=(1-u_1,0,\dots,0)^\top.
$$

Note that $M_n$ depends only on the transition probability vector
$p_n(\cdot)$, and hence $M_n\cdots M_1$ is the product of i.i.d.\
random (non-negative) matrices. By\break Furstenberg-Kesten's theorem,
the limiting behavior of such a product, as $n\to\infty$, is
controlled by the largest Lyapunov exponent
\begin{equation}\label{gamma1}
\gamma_1:=\lim_{n\to\infty}n^{-1}\ln \|M_n\dots M_1\|
\end{equation}
(by Kingman's sub-additive ergodic theorem, the\break limit exists
$\PP$-a.s.\ and is non-random). It follows that, $\Prob_0$-a.s.,
the RWRE $X_n$ is transient if and only if $\gamma_1\ne0$, and
moreover, $X_n\to+\infty\, (-\infty$) when $\gamma_1<0\,(\,>0)$,
whereas $\varliminf X_n=-\infty$, $\varlimsup X_n=+\infty$ when
$\gamma_1=0$.

For orientation, note that if $p_n(i)=p(i)$ are non-random
constants, then $\gamma_1=\ln\lambda_1$, where $\lambda_1>0$ is
the largest eigenvalue of $M_0$, and so $\gamma_1<0$ if and only
if $\lambda_1<1$. The latter means that the characteristic
polynomial $\varphi(\lambda):=\det (M_0-\lambda I)$ satisfies the
condition $(-1)^L\varphi(1)>0$. To evaluate $\det(M_0-I)$, replace
the first column by the sum of all columns and expand to get
$\varphi(1)=(-1)^{L-1} (b_1+\cdots +b_L)$. Substituting
expressions (\ref{eq:b}) it is easy to see that the above
condition amounts to $p(1)-\sum_{i=1}^L i\, p(-i)>0$, that is, the
mean drift of the random walk is positive and hence
$X_n\to+\infty$ a.s.

In the general case, $L\ge 1$, $R\ge 1$, similar considerations
lead to the following matrices of order $d:=L+R-1$ (cf.\ eqn
(\ref{eq:M}))
$$
\arraycolsep=.15pc
  M_n=\left(\begin{array}{cccccc}
 a_n(R-1)&\dots&a_n(1)& b_n(1)&\dots&b_n(L)\\[.3pc]
 1&0&\dots&\dots&\dots&0\\[.3pc]
 0&1&0&\dots&\dots&0\\[-.1pc]
 \vdots&\vdots&\ddots&\vdots&\vdots&\vdots\\[-.1pc]
 \vdots&\vdots&\vdots&\ddots&\vdots&\vdots\\[.1pc]
 0&\dots&\dots&0&1&0
 \end{array}\right)
$$
where $b_n(i)$ are given by (\ref{eq:b}) and
\begin{align*}
  {}&a_n(i):=-\frac{p_n(i)+\dots+p_n(R)}{p_n(R)}.
\end{align*}
Suppose that the ellipticity condition is satisfied in the form
$p_n(i)\ge\delta>0$, $i\ne0$, $-L\le i\le R$, and let
$\gamma_1\ge\gamma_2\ge\dots\ge\gamma_d$ be the (non-random)
Lyapunov exponents of $\{M_n\}$. The largest exponent $\gamma_1$
is again given by eqn (\ref{gamma1}), while other exponents are
determined recursively from the equalities
$$
\gamma_1+\dots+\gamma_k=\lim_{n\to\infty} n^{-1}
\ln\|{\wedge^k}(M_n\cdots M_1)\|
$$
($1\le k \le d$). Here $\wedge$ denotes the external
(anti-symmetric) product: $x\wedge y=-y\wedge x$
\,($x,y\in\RR^d$), and $\wedge^k M$ acts on the external product
space ${\wedge^k}\RR^d$, generated by the canonical basis
$\{e_{i_1}\wedge\cdots\wedge e_{i_k},\ 1\le i_1<\cdots<i_k\le
d\}$, as follows:
$$
{\wedge^k} M(x_1\wedge\cdots \wedge x_k):=M(x_1)\wedge\cdots
\wedge M(x_k).
$$

One can show that all exponents except $\gamma_R$ are sign
definite: $\gamma_{R-1}>0>\gamma_{R+1}$. Moreover, it is the sign
of $\gamma_R$ that determines whether the RWRE is transient or
recurrent, the dichotomy being the same as in the case $R=1$ above
(with $\gamma_1$ replaced by $\gamma_R$). Let us also mention that
an LLN and CLT can be proved here (see Br\'emont (2004)).

In conclusion, let us point out an alternative approach due to
Bolthausen and Goldsheid (2000) who studied a more general RWRE on a
strip $\ZZ\times \{0,1,\dots,\allowbreak m-1\}$. The link between
these two models is given by the representation $X_n=m Y_n+Z_n$,
where $m:=\max\{L,R\}$, \,$Y_n\in\ZZ$, \,$Z_n\in\{0,\dots,m-1\}$.
Random matrices arising here are constructed indirectly using an
auxiliary stationary sequence. Even though these matrices are
non-inde\-pen\-dent, thanks to their positivity the criterion of
transience can be given in terms of the sign of the \emph{largest
Lyapunov exponent}, which is usually much easier to deal with. An
additional attractive feature of this approach is that the condition
$p_x(R)>0$ ($\PP$-a.s.), which was essential for the previous
technique, can be replaced with a more natural condition
$\PP\{p_x(R)>0\}>0$.

\subsection{Random Bonds RWRE}\label{sec:Bonds}
Instead of having random probabilities of jumps at each \emph{site},
one could assign random weights to \emph{bonds} between the sites.
For instance, the transition probabilities $p_x=p(x,x+1,\omega)$ can
be defined by
\begin{equation}\label{a}
p_{x}=\frac{c_{x,x+1}}{c_{x-1,x}+c_{x,x+1}},
\end{equation}
where $c_{x,x+1}>0$ are i.i.d.\ random variables on the
environment space $\Omega$.

The difference between the two models may not seem very prominent,
but the behavior of the walk in the modified model (\ref{a})
appears to be quite different. Indeed, working as in Section
\ref{sec:TvR} we note that
$$
\rho_x=\frac{q_x}{p_x}=\frac{c_{x-1,x}}{c_{x,x+1}},
$$
hence, exploiting formulas (\ref{u1}) and (\ref{a}), we obtain,
$\PP$-a.s.,
\begin{equation}\label{eq:sym-rec}
\frac{1}{1-u_1}=\sum_{x=0}^{n-1} \frac{c_{01}}{c_{x,x+1}}\sim
c_{01}\, n\EE c_{01}^{-1}\to\infty,
\end{equation}
since $\EE c_{01}^{-1}>0$. Therefore, $f_{00}=1$, i.e.\ the random
walk is recurrent ($\Prob_0$-a.s.).

The method of environment viewed from the particle can also be
applied here (see Sznitman (2004)). Similarly to Section
\ref{sec:Lagrange}, we define a new probability measure
$\QQ=f(\omega)\PP$ using the density
$$
f(\omega)=Z^{-1}\bigl(c_{-1,0}(\omega)+c_{01}(\omega)\bigr),
$$
where $Z=2\EE c_{01}$ is the normalizing constant (we assume that
$\EE c_{01}<\infty$). One can check that $\QQ$ is invariant with
respect to the transition kernel eqn (\ref{a}), and by similar
arguments as in Section \ref{sec:Lagrange} we obtain that
$\lim_{n\to\infty}X_n/n$ exists ($\PW_0^\omega$-a.s.) and is given
by
$$
\int_{\Omega} d(0,\omega)\QQ(\dif\omega)
=Z^{-1}\EE\left[c_{01}-c_{-1,0}\right]=0,
$$
so the asymptotic velocity vanishes.

Furthermore, under suitable technical conditions on the
environment (e.g., $c_{01}$ being bounded away from $0$ and
$\infty$, cf.\ eqn (5)), one can prove the following CLT:
\begin{equation}\label{CLT-a}
  \lim_{n\to\infty}\Prob_0\biggl\{\frac{X_n}{\sqrt{n\sigma^2 }}\le
  x\biggr\}=\varPhi(x),
\end{equation}
where $\sigma^2=\left(\EE c_{01}\cdot \EE
c_{01}^{-1}\right)^{-1}$. Note that $\sigma^2\le 1$ (with a strict
inequality if $c_{01}$ is not reduced to a constant), which
indicates some slowdown in the spatial spread of the random bonds
RWRE, as compared to the ordinary symmetric random walk.

Thus, there is a dramatic distinction between the random bonds
RWRE, which is recurrent and diffusive, and the random sites RWRE,
with a much more complex asymptotics including both transient and
recurrent scenarios, slowdown effects and subdiffusive behavior.
This can be explained heuristically by noting that the random
bonds RWRE is reversible, that is, $m(x)\,p(x,y)=m(y)\,p(y,x)$ for
all $x,y\in\ZZ$, with $m(x):=c_{x-1,x}+c_{x,x+1}$ (this property
also easily extends to multidimensional versions). Hence, it
appears impossible to create extended traps which would retain the
particle for a very long time. Instead, the mechanism of the
diffusive slowdown in a reversible case is associated with the
natural variability of the environment resulting in the occasional
occurrence of isolated ``screening'' bonds with an anomalously
small weight $c_{x,x+1}$.

Let us point out that the RWRE determined by eqn (\ref{a}) can be
interpreted in terms of the random conductivity model (see Hughes,
1996). Suppose that each random variable $c_{x,x+1}$ attached to the
bond $(x,x+1)$ has the meaning of the conductance of this bond (the
reciprocal, $c_{x,x+1}^{-1}$, being its resistance). If a voltage
drop $V$ is applied across the system of $N$ successive bonds, say
from $0$ to $N$, then the same current $I$ flows in each of the
conductors and by\break Ohm's law we have $I=c_{x,x+1} V_{x,x+1}$,
where $V_{x,x+1}$ is the voltage drop across the corresponding bond.
Hence
$$
  V=\sum_{x=0}^N V_{x,x+1}=I\sum_{x=0}^N c_{x,x+1}^{-1},
$$
which amounts to saying that the total resistance of the system of
consecutive elements is given by the sum of the individual
resistances. The effective conductivity of the finite system,
$\overline{c}_N$, is defined as the average conductance per bond,
so that
$$
  \overline{c}^{\,-1}_N=
  \frac{1}{N}\sum_{x=0}^N c_{x,x+1}^{-1},
$$
and by the strong LLN, $\overline{c}_N^{\,-1}\to\EE c_{01}^{-1}$
as $N\to\infty$ ($\PP$-a.s.). Therefore, the effective
conductivity of the infinite system is given by
$\overline{c}=\left(\EE c_{01}^{-1}\right)^{-1}$, and we note that
$\overline{c}<\EE c_{01}$ if the random medium is non-degenerate.

Returning to the random bonds RWRE, eqn (\ref{a}), it is easy to
see that a site $j$ is recurrent if and only if the conductance
$c_{j,\infty}$ between $x$ and $\infty$ equals zero. Using again
Ohm's law, we have (cf.\ eqn (\ref{eq:sym-rec}))
$$
c^{-1}_{j,+\infty}=\sum_{x=j}^\infty c_{x,x+1}^{-1}= \infty,\quad
\PP\text{-a.s.}
$$
and we recover the result about recurrence.

\subsection{Continuous-Time RWRE}\label{sec:continuous}
As in the discrete-time case, a random walk on $\ZZ$ with continuous time
is a homogeneous Markov chain $X_t$, $t\in[0,\infty)$, with state space $\ZZ$
and nearest neighbor (or at least bounded) jumps. The term ``Markov'' as usual
refers to the ``lack of memory'' property, which amounts to saying that
from the entire history of the process development up to a given time, only
the current position of the walk is important for the future evolution while
all other information is irrelevant.

Since there is no smallest time unit as in the dis\-crete-time
case, it is convenient to describe transitions of $X_t$ in terms
of transition rates characterizing the likelihood of various jumps
during a very short time. More precisely, if
$p_{xy}(t):=P\{X_t=y\,|\,X_0=x\}$ are the transition probabilities
over time $t$, then for $h\to0$
\begin{equation}\label{a1}
\begin{aligned}
p_{xy}(h)&=c_{xy} h + o(h)\qquad (x\ne y),\\
p_{xx}(h)&=1-h\sum_{y\ne x} c_{xy} + o(h).
\end{aligned}
\end{equation}
Equations for the functions $p_{xy}(t)$ can then be derived by
adapting the method of decomposition commonly used for
dis\-crete-time Mar\-kov chains (cf.\ Section \ref{sec:TvR}). Here
it is more convenient to decompose with respect to the ``last''
step, i.e.\ by considering all possible transitions during a small
increment of time at the end of the time interval $[0,t+h]$. Using
Markov property and eqn (\ref{a1}) we can write
\begin{align*}
p_{0x}(t+h)
&=h\sum_{y\ne x} p_{0y}(t)\,c_{yx}\\
&+ p_{0x}(t) \biggl(1-h\sum_{y\ne x} c_{xy}\biggr) +o(h),
\end{align*}
which in the limit $h\to0$ yields the
\emph{master equation} (or \emph{Chapman-Kolmogorov's forward
equation})
\begin{equation}\label{master}
\begin{aligned}
{}&\frac{\dif}{\dif t}p_{0x}(t)=\sum_{y\ne x}
\bigl\{c_{yx}p_{0y}(t)-c_{xy}p_{0x}(t)\bigr\} ,\\
&p_{0x}(0)=\delta_0(x),
\end{aligned}
\end{equation}
where $\delta_0(x)$ is the Kronecker symbol.

Continuous-time RWRE are therefore naturally described via the
randomized master equation, i.e.\ with random transition rates.
The canonical example, originally motivated by Dyson's study of
the chain of harmonic oscillators with random couplings, is a
symmetric nearest-neighbor RWRE, where the random transition rates
$c_{xy}$ are non-zero only for $y=x\pm1$ and satisfy the condition
$c_{x,x+1}=c_{x+1,x}$, otherwise being i.i.d.\ (see Alexander
\emph{et al.} (1981)). In this case, the problem (\ref{master}) can
be formally solved using the Laplace transform, leading to the
equations
\begin{align}
\label{hat_p1}
s+G^{+}_0+G^{-}_0&=[\hat p_{0}(s)]^{-1},\\
\label{hat_p2} s+G^{-}_x + G^{+}_x&=0\quad (x\ne 0),
\end{align}
where $G^{-}_x$, $G^{+}_x$ are defined as
\begin{equation}\label{G_pm}
G^{\pm}_x:=c_{x,x\pm1}\, \frac{\hat p_{0x}(s)-\hat
p_{0,x\pm1}(s)}{\hat p_{0x}(s)}
\end{equation}
and $\hat p_{0x}(s):=\int_0^\infty p_{0x}(t)\,\e^{-st}\,\dif t$.
From eqs (\ref{hat_p2}), (\ref{G_pm}) one obtains the recursion
\begin{align}
\label{G+}
G^{\pm}_x&=\biggl(\frac{1}{c_{x,x\pm1}}+\frac{1}{s+G_{x\pm1}^{+}}\biggr)^{-1},\\
\notag
x&=0,\pm1,\pm2,\dots
\end{align}
The quantities $G^{\pm}_0$ are therefore expressed as infinite
continued fractions depending on $s$ and the random variables
$c_{x,x\pm1}$, $c_{x,x\pm2}$, \dots The function $\hat p_{00}(s)$
can then be found from eqn (\ref{hat_p1}).

In its generality, the problem is far too hard, and we shall only comment on
how one can evaluate the annealed mean
$$
\EE \hat p_{00}(s)=\EE (s+G^{+}_0+G^{-}_0)^{-1}.
$$
According to eqn (\ref{G+}), the random variables $G^{+}_0$,
$G^{-}_0$ are determined by the same algebraic formula, but
involve the rate coefficients from different sides of site $x$,
and hence are i.i.d. Furthermore, eqn (\ref{G+}) implies that the
random variables $G^{+}_0$, $G^{+}_{1}$ have the same distribution
and, moreover, $G^{+}_{1}$ and $c_{01}$ are independent.
Therefore, eqn (\ref{G+}) may be used as an integral equation for
the unknown density function of $G^{+}_{0}$. It can be proved that
the suitable solution exists and is unique, and although an
explicit solution is not available, one can obtain the asymptotics
of small values of $s$, thereby rendering information about the
behavior of $p_{00}(t)$ for large $t$. More specifically, one can
show that if $c_*:=\bigl(\EE c_{01}^{-1}\bigr)^{-1}>0$ then
$$
\EE \hat p_{00}(s)\sim (4c_* s)^{-1/2}, \quad s\to0,
$$
and so by a Tauberian theorem
\begin{equation}\label{c_*}
 \EE p_{00}(t)\sim (4\pi c_* t)^{-1/2},
 \quad t\to\infty.
\end{equation}
Note that asymptotics (\ref{c_*}) appears to be the same as for an
ordinary symmetric random walk with constant transition rates
$c_{x,x+1}=c_{x+1,x}= c_*$, suggesting that the latter provides an
``effective medium approximation'' (EMA) for the RWRE considered
above.

This is further confirmed by the asymptotic calculation of the
annealed mean-square displacement, $\Exp_0 X_t^2 \sim 2c_* t$ as
$t\to\infty$ (Alexander \emph{et al.}\ 1981). Moreover, Kawazu and
Kesten (1984) proved that $X_t$ is asymptotically normal:
\begin{equation}\label{CLT-a1}
  \lim_{t\to\infty}\Prob_0 \biggl\{\frac{X_t}{\sqrt{2c_*t}}\le x\biggr\}=
   \varPhi(x).
\end{equation}
Therefore, if $c_*>0$ then the RWRE has the same diffusive
behavior as the corresponding ordered system, with a well-defined
diffusion constant $D=c_*$.

In the case where $c_*=0$ (i.e., $\EE c_{01}^{-1}=\infty$), one
may expect that the RWRE exhibits subdiffusive behavior. For
example, if the density function of the transition rates is
modelled by
$$
f_a(u)=(1-\alpha)\, u^{-\alpha}{\mathbf 1}_{\{0<u<1\}}\quad
(0<\alpha<1),
$$
then, as shown by Alexander \emph{et al.} (1981),
\begin{align*}
  \EE p_{00}(t)&\sim C_\alpha\, t^{-(1-\alpha)/(2-\alpha)},\\
  \Exp_0 X_t^2&\sim C'_\alpha\, t^{2(1-\alpha)/(2-\alpha)}.
\end{align*}
In fact, Kawazu and Kesten (1984) proved that in this case
$t^{-\alpha/(1+\alpha)}X_t$ has a (non-Gaussian) limit distribution
as $t\to\infty$.

To conclude the discussion of the continuous-time case, let us
point out that some useful information about recurrence of $X_t$
can be obtained by considering an \emph{imbedded} (discrete-time)
random walk $\tilde X_n$, defined as the position of $X_t$ after
$n$ jumps. Note that continuous-time Markov chains admit an
alternative description of their evolution in terms of sojourn
times and the distribution of transitions at a jump. Namely, if
the environment $\omega$ is fixed then the random sojourn time of
$X_t$ in each state $x$ is exponentially distributed with mean
$1/c_x$, where $c_x:=\sum_{y\ne x}c_{xy}$, while the distribution
of transitions from $x$ is given by the probabilities
$p_{xy}=c_{xy}/c_x$.

For the symmetric nearest-neighbor RWRE considered above, the
transition probabilities of the imbedded random walk are given by
\begin{align*}
  p_x:=p_{x,x+1}&=\frac{c_{x,x+1}}{c_{x-1,x}+c_{x,x+1}},\\
  q_x:=p_{x,x-1}&=1-p_x,
\end{align*}
and we recognize here the transition law of a random walk in the
random bonds environment considered in Section \ref{sec:Bonds}
(cf.\ eqn (\ref{a})). Recurrence and zero asymptotic velocity
established there are consistent with the results discussed in the
present section (e.g., note that the CLT for both $X_n$, eqn
(\ref{CLT-a}), and $X_t$, eqn (\ref{CLT-a1}), does not involve any
centering). Let us point out, however, that a ``naive''
discretization of time using the mean sojourn time appears to be
incorrect, as this would lead to the scaling $t=n\delta_1$ with
$\delta_1:=\EE (c_{-1,0}+c_{01})^{-1}$, while from comparing the
limit theorems in these two cases, one can conclude that the true
value of the effective discretization step is given by
$\delta_*:=(2c_*)^{-1}=\frac12 \EE c_{01}^{-1}$. In fact, by the
arith\-metic-har\-monic mean inequality we have
$\delta_*>\delta_1$, which is a manifestation of the RWRE's
diffusive slowdown.

\section{RWRE in Higher Dimensions}
Multidimensional RWRE with nearest-neighbor\break jumps are defined in a
similar fashion: from site $x\in\ZZ^d$ the random walk can jump to
one of the $2d$ adjacent sites $x+e\in\ZZ^d$ (such that $|e|=1$),
with probabilities $p_x(e)\ge0$, \,$\sum_{|e|=1}p_x(e)=1$, where
the random vectors $p_x(\cdot)$ are assumed to be i.i.d.\ for
different $x\in\ZZ^d$. As usual, we will also impose the condition
of uniform ellipticity:
\begin{equation}\label{eq:ellipt}
p_x(e)\ge\delta>0,\quad |e|=1, \ \ x\in\ZZ^d,\ \,\PP\text{-a.s.}.
\end{equation}

In contrast to the one-dimensional case, theory of RWRE in higher
dimensions is far from maturity. Possible asymptotic behaviors of
the RWRE for $d\ge 2$ are not understood well enough, and many basic
questions remain open. For instance, no definitive classification of
the RWRE is available regarding transience and recurrence.
Similarly, LLN and CLT have been proved only for a limited number of
specific models, while no general sharp results have been obtained.
On a more positive note, there has been considerable progress in
recent years in the so-called ballistic case, where powerful
techniques have been developed (see Sznitman (2002, 2004) and
Zeitouni (2003, 2004)). Unfortunately, not much is known for
non-ballistic RWRE, apart from special cases of balanced RWRE in
$d\ge 2$ (Lawler 1982), small isotropic perturbations of ordinary
symmetric random walks in $d\ge 3$ (Bricmont and Kupiainen 1991),
and some examples based on combining components of ordinary random
walks and RWRE in $d\ge7$ (Bolthausen \emph{et al.}\ 2003). In
particular, there are no examples of subdiffusive behavior in any
dimension $d\ge 2$, and in fact it is largely believed that a CLT is
always true in any uniformly elliptic, i.i.d.\ random environment in
dimensions $d\ge3$, with somewhat less certainty about $d=2$. A
heuristic explanation for such a striking difference with the case
$d=1$ is that due to a less restricted topology of space in higher
dimensions, it is much harder to force the random walk to visit
traps, and hence the slowdown is not so pronounced.

In what follows, we give a brief account of some of the known
results and methods in this fast developing area (for further
information and specific references, see an extensive review by
Zeitouni (2004)).

\subsection{Zero-One Laws and LLNs}\label{sec:0-1}
A natural first step in a multidimensional context is to explore
the behavior of the random walk $X_n$ as projected on various
one-dimensional straight lines. Let us fix a test unit vector
$\ell\in\RR^d$, and consider the process
$Z_n^{\ell}:=X_n\cdot\ell$. Then for the events
$A_{\pm\ell}:=\{\lim_{n\to\infty} Z^\ell_n=\pm\infty\} $ one can
show that
\begin{equation}\label{eq:0-1}
\Prob_0(A_\ell\cup A_{-\ell})\in\{0,1\}.
\end{equation}
That
is to say, for each $\ell$ the probability that the random walk
escapes to infinity in the direction $\ell$ is either $0$ or $1$.

Let us sketch the proof. We say that $\tau$ is \emph{record time}
if $|Z^\ell_t|>|Z^\ell_k|$ for all $k<t$, and \emph{regeneration
time} if in addition $|Z^\ell_{\tau}|\le |Z^\ell_n|$ for all $n\ge
\tau$. Note that by the ellipticity condition (\ref{eq:ellipt}),
$\varlimsup_{n\to\infty} |Z^\ell_n|= \infty$ ($\Prob_0$-a.s.),
hence there is an infinite sequence of record times
$0=\tau_0<\tau_1<\tau_2<\cdots$ If $\Prob_0(A_\ell\cup
A_{-\ell})>0$, we can pick a subsequence of record times
$\tau'_i$, each of which has a \emph{positive}
$\Prob_0$-probability to be a regeneration time (because otherwise
$|Z^\ell_n|$ would persistently backtrack towards the origin and
the event $A_\ell\cup A_{-\ell}$ could not occur). Since the
trials for different record times are independent, it follows that
a regeneration time $\tau^*$ occurs $\Prob_0$-a.s. Repeating this
argument, we conclude that there exists an infinite sequence of
regeneration times $\tau_i^*$, which implies that
$|Z^\ell_n|\to\infty$ ($\Prob_0$-a.s.), i.e., ${\Prob(A_\ell\cup
A_{-\ell})=1}$.

Regeneration structure introduced by the sequence $\{\tau_i^*\}$
plays a key role in further analysis of the RWRE and is
particularly useful for proving an LLN and a CLT, due to the fact
that pieces of the random walk between consecutive regeneration
times (and fragments of the random environment involved thereby)
are independent and identically distributed (at least starting
from $\tau_1^*$). In this vein, one can prove a ``directional''
version of the LLN, stating that for each $\ell$ there exist
deterministic $v_\ell,v_{-\ell}$ (possibly zero) such that
\begin{equation}\label{eq:Z/n}
  \lim_{n\to\infty}\frac{Z_n^\ell}{n}=v_\ell\,{\mathbf 1}_{A_\ell}+
  v_{-\ell}\,{\mathbf 1}_{A_{-\ell}},\quad \Prob_0\text{-a.s.}
\end{equation}

Note that if $\Prob_0(A_\ell)\in\{0,1\}$, eqn (\ref{eq:Z/n}) in
conjunction with eqn (\ref{eq:0-1}) would readily imply
\begin{equation}\label{eq:Z/n1}
  \lim_{n\to\infty}\frac{Z_n^\ell}{n}=v_\ell, \quad \Prob_0\text{-a.s.}
\end{equation}
Moreover, if $\Prob_0(A_{\ell})\in\{0,1\}$ for any $\ell$, then
there exists a deterministic $v$ (possibly zero) such that
\begin{equation}\label{eq:Z/n2}
  \lim_{n\to\infty}\frac{X_n}{n}=v, \quad \Prob_0\text{-a.s.}
\end{equation}
Therefore, it is natural to ask if a zero-one law (\ref{eq:0-1})
can be enhanced to that for the individual probabilities
$\Prob_0(A_\ell)$. It is known that the answer is affirmative for
i.i.d.\ environments in $d=2$, where indeed
$\Prob(A_\ell)\in\{0,1\}$ for any $\ell$, with counter-examples in
certain stationary ergodic (but not uniformly elliptic)
environments. However, in the case $d\ge3$ this is an open
problem.

\subsection{Kalikow's Condition and Sznitman's\\ Condition
(\boldmath{$\mathbf T'$})}\label{sec:Kalikow} An RWRE is called
\emph{ballistic} (\emph{ballistic in direction $\ell$}) if $v\ne0$
($v_\ell\ne0$), see eqs (\ref{eq:Z/n1}), (\ref{eq:Z/n2}). In this
section, we describe conditions on the random environment which
ensure that the RWRE is ballistic.

Let $U$ be a connected strict subset of $\ZZ^d$ containing the
origin. For $x\in U$, denote by
$$
g(x,\omega):=\EW^\omega_0 \sum_{n=0}^{T_U}
\mathbf 1_{\{X_n=x\}}
$$
the quenched mean number of visits to $x$ prior to the exit time
$T_U:=\min\{n\ge0: X_n\notin U\}$. Consider an auxiliary Markov
chain $\widehat X_n$, which starts from $0$, makes
nearest-neighbor jumps while in $U$, with (non-random)
probabilities
\begin{equation}\label{eq:Xhat}
\begin{aligned}
\displaystyle
\widehat p_x(e)&=
 \frac{\EE \left[g(x,\omega) p_x(e)\right]}
      {\EE \left[g(x,\omega)\right]},\quad
x\in U,
\end{aligned}
\end{equation}
and is absorbed as soon as it first leaves $U$. Note that the
expectations in eqn (\ref{eq:Xhat}) are finite; indeed, if
$\alpha_x$ is the probability to return to $x$ before leaving $U$,
then, by the Markov property, the mean number of returns is given
by
$$
\sum_{k=1}^\infty k \alpha_x^k(1-\alpha_x)=\frac{\alpha_x}{1-\alpha_x}<\infty,
$$
since, due to ellipticity, $\alpha_x<1$.

An important property, highlighting the usefulness of $\widehat
X_n$, is that if $\widehat X_n$ leaves $U$ with probability $1$,
then the same is true for the original RWRE $X_n$ (under the
annealed law $\Prob_0$), and moreover, the exit points $\widehat
X_{{\widehat T}_U}$ and $X_{T_U}$ have the same distribution laws.

Let $\ell\in\RR^d$, $|\ell|=1$. One says that \emph{Kalikow's
condition with respect to $\ell$} holds if the local drift of
$\widehat X_n$ in the direction $\ell$ is uniformly bounded away
from zero:
\begin{equation}\label{eq:Kalikow}
  \inf_U\inf_{x\in U}
  \sum_{|e|=1} (e\cdot \ell)\,\widehat p_x(e)>0.
\end{equation}
A sufficient condition for (\ref{eq:Kalikow}) is, for example,
that for some $\kappa>0$
$$
  \EE\left[(d(0,\omega)\cdot\ell)_+\right]\ge \kappa
  \EE\left[(d(0,\omega)\cdot\ell)_-\right],
$$
where $d(0,\omega)=\EW_0^\omega X_1$ and $u_\pm:=\max\{\pm u,0\}$.

A natural implication of Kalikow's condition (\ref{eq:Kalikow}) is that
$\Prob_0(A_\ell)=1$ and $v_\ell>0$ (see eqn (\ref{eq:Z/n1})). Moreover, noting
that eqn (\ref{eq:Kalikow}) also holds for all $\ell'$ in a vicinity of $\ell$
and applying the above result with $d$ non-collinear vectors from that
vicinity, we conclude that under Kalikow's condition there exists a
deterministic $v\ne0$ such that $X_n/n\to v$ as $n\to\infty$ ($\Prob_0$-a.s.).
Furthermore, it can be proved that $(X_n-nv)/\sqrt{n}$ converges in law to a
Gaussian distribution (see Sznitman (2004)).

It is not hard to check that in dimension $d=1$ Kalikow's
condition is equivalent to $v\ne0$ and therefore characterizes
completely all ballistic walks. For $d\ge 2$, the situation is
less clear; for instance, it is not known if there exist RWRE with
$\Prob(A_\ell)>0$ and $v_\ell=0$ (of course, such RWRE cannot
satisfy Kalikow's condition).

Sznitman (2004) has proposed a more complicated transience
condition ($\text{T}'$) involving certain regeneration times
$\tau_i^*$ similar to those described in Section \ref{sec:0-1}. An
RWRE is said to satisfy \emph{Sznitman's condition} ($\text{T}'$)
\emph{relative to direction $\ell$} if $\Prob_0(A_\ell)=1$ and for
some $c>0$ and all $0<\gamma<1$
\begin{equation}\label{eq:T'}
  \Exp_0 \exp\Bigl(c\sup_{n\le\tau_1^*} |X_n|^\gamma\Bigr)<\infty.
\end{equation}
This condition provides a powerful control over $\tau_1^*$ for
$d\ge 2$ and in particular ensures that $\tau_1^*$ has finite
moments of any order. This is in sharp contrast with the
one-dimensional case, and should be viewed as a reflection of much
weaker traps in dimensions $d\ge 2$. Condition (\ref{eq:T'}) can
also be reformulated in terms of the exit distribution of the RWRE
from infinite thick slabs ``orthonormal'' to directions $\ell'$
sufficiently close to $\ell$. As it stands, the latter
reformulation is difficult to check, but Sznitman (2004) has
developed a remarkable ``effective'' criterion reducing the job to
a similar condition in \emph{finite} boxes, which is much more
tractable and can be checked in a number of cases.

In fact, condition ($\text{T}'$) follows from Kalikow's condition,
but not the other way around. In the one-dimensional case,
condition ($\text{T}'$) (applied to $\ell=1$ and $\ell=-1$) proves
to be equivalent to the transient behavior of the RWRE, which, as
we have seen in Theorem \ref{Solomon2} (Section \ref{sec:Stokes}),
may happen with $v=0$, i.e.\ in a non-ballistic scenario. The
situation in $d\ge 2$ is quite different, as condition
($\text{T}'$) implies that the RWRE is ballistic in the direction
$\ell$ (with $v_\ell>0$) and satisfies a CLT (under $\Prob_0$). It
is not known whether the ballistic behavior for $d\ge 2$ is
completely characterized by condition ($\text{T}'$), although this
is expected to be true.

\subsection{Balanced RWRE}\label{sec:Balanced}
In this section we discuss a particular case of non-ballistic
RWRE, for which LLN and CLT can be proved. Following Lawler
(1982), we say that an RWRE is \emph{balanced} if $p_x(e)=p_x(-e)$
for all $x\in\ZZ^d$, $|e|=1$ ($\PP$-a.s.). In this case, the local
drift vanishes, $d(x,\omega)=0$, hence the coordinate processes
$X_n^i$ ($i=1,\dots,d$) are martingales with respect to the
natural filtration $\calF_n=\sigma\{X_0,\dots,X_n\}$. The quenched
covariance matrix of the increments $\Delta
X^i_n:=X_{n+1}^i-X_{n}^i$ ($i=1,\dots,d$) is given by
\begin{equation}\label{eq:Delta}
\EW_0^\omega\left[\Delta X_n^i\, \Delta X_n^j\,|\,\calF_{n}\right]
=2\delta_{ij} p_{X_{n}}(e_i).
\end{equation}
Since the right-hand side of eqn (\ref{eq:Delta}) is uniformly
bounded, it follows that $X_n/n\to 0$ ($\Prob_0$-a.s.). Further,
it can be proved that there exist deterministic positive constants
$a_1,\dots,a_d$ such that for $i=1,\dots,d$
\begin{equation}\label{eq:a}
 \lim_{n\to\infty}\frac{1}{n}\sum_{k=0}^{n-1} p_{X_k}(e_i)=\frac{a_i}{2},\quad
\Prob_0\text{-a.s.}
\end{equation}
Once this is proved, a multidimensional CLT for martingale
differences yields that $X_n/\sqrt{n}$ converges in law to a
Gaussian distribution with zero mean and the covariances
$b_{ij}=\delta_{ij}a_i$.

The proof of (\ref{eq:a}) employs the method of environment viewed
from the particle (cf.\ Section \ref{sec:Lagrange}). Namely,
define a Markov chain $\omega_n:=\theta^{X_n}\omega$ with the
transition kernel
\begin{align*}
T(\omega,\dif\omega')=\sum_{i=1}^d
\bigl[&p_0(e_i)\,\delta_{\theta\omega}(\dif\omega')\\
&+p_0(-e_i)\,\delta_{\theta^{-1}\omega}(\dif\omega')\bigr]
\end{align*}
(cf.\ eqn (\ref{eq:T})). The next step is to find a probability
measure $\QQ$ on $\Omega$ invariant under $T$ and absolutely
continuous with respect to $\PP$. Unlike the one-dimensional case,
however, an explicit form of $\QQ$ is not available, and $\QQ$ is
constructed indirectly as the limit of invariant measures of
certain periodic modifications of the RWRE. Birkhoff's ergodic
theorem then yields, $\Prob_0$-a.s.,
\begin{align*}
\frac{1}{n}\sum_{k=0}^{n-1} p_{X_k}(e_i,\omega)&=
\frac{1}{n}\sum_{k=0}^{n-1} p_0(e_i,\omega_k)\\
&\to\int_\Omega p_0(e_i,\omega_1)\QQ(\dif\omega)\ge\delta,
\end{align*}
by the ellipticity condition (\ref{eq:ellipt}), and eqn
(\ref{eq:a}) follows.

With regard to transience, balanced RWRE admit a complete and simple
classification. Namely, it has been proved (see Zeitouni (2004)) that any
balanced RWRE is transient for $d\ge3$ and recurrent for $d=2$
($\Prob_0$-a.s.). It is interesting to note, however, that these answers may be
false for certain balanced random walks in a \emph{fixed} environment
($\PP$-probability of such environments being zero, of course). Indeed, examples
can be constructed of balanced random walks in $\ZZ^2$ and in $\ZZ^d$ with
$d\ge 3$, which are transient and recurrent, respectively (Zeitouni 2004).

\subsection{RWRE Based on Modification of\\ Ordinary Random Walks}
A number of partial results are known for RWRE constructed on the
basis of ordinary random walks via certain randomization of the
environment. A natural model is obtained by a small perturbation
of a simple symmetric random walk. To be more precise, suppose
that: (a) $|p_x(e)-\frac{1}{2d}|<\varepsilon$ for all $x\in\ZZ^d$
and any $|e|=1$, where $\varepsilon>0$ is small enough; (b) $\EE
p_x(e)=\frac{1}{2d}\strut$; (c) vectors $p_x(\cdot)$ are i.i.d.\
for different $x\in\ZZ^d$,
and (d) the distribution of the vector $p_x(\cdot)$ is isotropic,
i.e.\ invariant with respect to permutations of its coordinates.
Then for $d\ge 3$ Bricmont and Kupiainen (1991) have proved an LLN
(with zero asymptotic velocity) and a quenched CLT (with
non-degenerate covariance matrix). The proof is based on the
renormalization group method, which involves decimation in time
combined with a suitable spatial-temporal scaling. This
transformation replaces an RWRE by another RWRE with weaker
randomness, and it can be shown that iterations converge to a
Gaussian fixed point.

Another class of examples are also built using small perturbations
of simple symmetric random walks, but are \emph{anisotropic} and
exhibit ballistic behavior, providing that the annealed local
drift in some direction is strong enough (see Sznitman (2004)).
More precisely, suppose that $d\ge 3$ and $\eta\in(0,1)$. Then
there exists $\varepsilon_0=\varepsilon_0(d,\eta)>0$ such that if
$|p_x(e)-\frac{1}{2d}\strut|<\varepsilon$ \,($x\in\ZZ^d$, $|e|=1$)
with $0<\varepsilon<\varepsilon_0$, and for some $e_0$ one has
$\EE{}[d(x,\omega)\cdot e_0]\ge \varepsilon^{2.5-\eta}$ ($d=3$) or
$\ge \varepsilon^{3-\eta}$ ($d\ge 4$), then Sznitman's condition
($\text{T}'$) is satisfied with respect to $e_0$ and therefore the
RWRE is ballistic in the direction $e_0$ (cf.\ Section
\ref{sec:Kalikow}).

Examples of a different type are constructed in dimensions $d\ge
6$ by letting the first $d_1\ge 5$ coordinates of the RWRE $X_n$
behave according to an ordinary random walk, while the remaining
$d_2=d-d_1$ coordinates are exposed to a random environment (see
Bolthausen \emph{et al.} (2003)). One can show that there exists a
deterministic $v$ (possibly zero) such that $X_n/n\to v$
($\Prob_0$-a.s.). Moreover, if $d_1\ge 13$ then
$(X_n-nv)/\sqrt{n}$ satisfies both quenched and annealed CLT.
Incidentally, such models can be used to demonstrate the
surprising features of the multidimensional RWRE. For instance,
for $d\ge 7$ one can construct an RWRE $X_n$ such that the
annealed local drift does not vanish, $\EE d(x,\omega)\ne0$, but
the asymptotic velocity is zero, $X_n/n\to0$ ($\Prob_0$-a.s.), and
furthermore, if $d\ge15$ then in this example $X_n/\sqrt{n}$
satisfies a quenched CLT. (In fact, one can construct such RWRE as
small perturbations of a simple symmetric walk.) On the other
hand, there exist examples (in high enough dimensions) where the
walk is ballistic with a velocity which has an opposite direction
to the annealed drift $\EE d(x,\omega)\ne0$. These striking
examples provide ``experimental'' evidence of many unusual
properties of the multidimensional RWRE, which, no doubt, will be
discovered in the years to come.

\bigskip\noindent
{{\it See also:} \,Averaging Methods; Growth Processes in Random
Matrix Theory; Lagrangian Dispersion\break (Passive Scalar); Random
Dynamical Systems; Random Matrix Theory in Physics; Stochastic
Differential Equations; Stochastic Loewner Evolutions.

}

\end{document}